# GROUPES DE SELMER ET ACCOUPLEMENTS ; CAS PARTICULIER DES COURBES ELLIPTIQUES

BERNADETTE PERRIN-RIOU

Soit $E$ une courbe elliptique définie sur le corps des nombres rationnels $\mathbb{Q}$. D'après le théorème de Mordell, le groupe $E(\mathbb{Q})$ des points de $E$ rationnels sur $\mathbb{Q}$ est un $\mathbb{Z}$-module de type fini. Nekovář vient de démontrer que son rang est de même parité que la multiplicité du zéro en $s = 1$ de la fonction $L$ complexe associée à $E/\mathbb{Q}$ lorsque le groupe de Tate-Shafaravich est fini. La conjecture de Birch et Swinnerton-Dyer prédit qu'il y a égalité entre ces deux entiers attachés à $E$.

La démonstration de ce résultat par Nekovář utilise essentiellement trois arguments et se fait en introduisant un corps quadratique imaginaire $K$ et un nombre premier $p$ auxiliaires vérifiant certaines conditions. Le premier argument utilise le théorème de Cornut concernant les points de Heegner ([8], conjecture de Mazur [20]) : on en déduit que le rang de $E(K_n)$ tend vers l'infini avec $n$ pour $K_n$ l'extension de $\mathbb{Q}$ de de degré $2p^n$ et diédrale sur $\mathbb{Q}$. Le deuxième argument est à base de théorie d'Iwasawa. Il s'agit de généraliser le théorème de Cassels qui affirme que le groupe de Tate-Shafaravich d'une courbe elliptique modulo sa partie divisible (noté div) est un carré, ou ce qui revient presque au même, que le quotient du $p$-groupe de Selmer $S_p(E/K)$ de $E/K$ modulo sa partie divisible est un carré : on peut construire une forme alternée et non dégénérée sur $S(K)/\operatorname{div}$. Ce qui joue ici le rôle du groupe de Tate-Shafarevich est le quotient de $S_p(E/K_\infty)$ par sa partie $\Lambda$-divisible $\operatorname{div}_\Lambda$ où $K_\infty = \cup K_n$, $\Gamma = \operatorname{Gal}(K_\infty/K)$ et $\Lambda$ est l'algèbre de groupe continue $\mathbb{Z}_p[[\Gamma]]$. On peut encore construire sur $S_p(E/K_\infty)/\operatorname{div}_\Lambda$ une forme $\Lambda$-linéaire et alternée. Le troisième argument utilise les résultats de Kolyvagin généralisés par Bertolini et Darmon ([17], [4]) et des arguments de descente pour conclure.

Pour construire la forme alternée, Nekovář reprend complètement la théorie des groupes de Selmer en utilisant le formalisme des complexes. Il obtient ainsi d'autres applications en théorie de Hida et autres. Nous nous contentons ici de faire cette construction en allant au plus court et de replacer ensuite ces résultats dans un contexte plus général.

Le principe de la démonstration est de faire grand usage du twist d'un $\Lambda$-module : l'adjoint d'un $\Lambda$-module de torsion se calcule en effet facilement lorsque ses coinvariants sont de torsion pour l'anneau quotient, ce qui est réalisable en faisant un twist convenable. Cette astuce permet d'éviter les difficultés dues au fait que le groupe de Mordell-Weil n'est pas fini. Ainsi, par exemple, l'accouplement de Cassels-Tate peut se calculer comme une "limite convenable" des accouplements relatifs aux représentations twistées par $k$.

## Table des matières

1. Préliminaires d'algèbre commutative

1.1. **Adjoint et dualité.** Soit $\Lambda$ un anneau local noetherien complet de dimension $r$ ; plus précisément $\Lambda$ sera l'algèbre de groupes continue d'un groupe $\Gamma$ topologiquement isomorphe à $\mathbb{Z}_p^{r-1}$ : $\Lambda = \mathbb{Z}_p[[\Gamma]]$. Soit $X$ un $\Lambda$-module de type fini. On pose pour tout entier $i \geq 0$

$$a_\Lambda^i(M) = \text{Ext}_\Lambda^i(M, \Lambda)$$

En particulier le $\Lambda$-module $a_\Lambda^1(M) = \text{Ext}_\Lambda^1(M, \Lambda)$ est l'adjoint (d'Iwasawa) de $M$. Rappelons quelques faits d'algèbre commutative.
- Un module est dit pseudo-nul si la hauteur des idéaux associés est supérieure ou égale à 2 ;
- La hauteur des idéaux associés à $a_\Lambda^i(M)$ est supérieure ou égale à $i$ ; aussi, $a_\Lambda^i(M)$ est pseudo-nul pour $i \geq 2$, $a_\Lambda^1(M)$ est un $\Lambda$-module de torsion et pour $\dim \Lambda = 3$, $a_\Lambda^3(M)$ est fini ;
- On peut interpréter $a_\Lambda^1(M)$ de la manière suivante. Soit $\text{Frac}\,\Lambda$ le corps des fractions de $\Lambda$. Alors,

$$a_\Lambda^1(M) = \text{Hom}_\Lambda(M, \text{Frac}\,\Lambda/\Lambda) \ .$$

- Si $M$ est un $\Lambda$-module, $a_\Lambda^1(M)$ n'a pas de sous-modules pseudo-nuls non nuls. En particulier, si $M$ est un module pseudo-nul, $a_\Lambda^1(M) = 0$.
- Supposons $\Lambda$ de dimension 3. Si $M$ n'a pas de sous-modules pseudonuls non finis, alors $a_\Lambda^2(M)$ est fini.

Démontrons les deux derniers points. Il existe un homomorphisme injectif $E \to M$ à conoyau $F$ pseudonul avec $E$ un module élémentaire $\oplus \Lambda/(f_i)$. On en déduit la suite exacte :

$$0 \to a_\Lambda^1(M) \to a_\Lambda^1(E) \to a_\Lambda^2(F) \ .$$

Il est facile de voir que $E$ est de dimension projective 1, ce qui implique que $a_\Lambda^1(M)$ n'a pas de sous-modules pseudo-nuls non nuls. D'autre part, il existe un homomorphisme $M \to E$ à conoyau pseudo-nul $F$ avec $E$ module élémentaire et de noyau pseudo-nul. Supposons maintenant que la dimension projective de $\Lambda$ est 3 et que $M$ n'a pas de sous-modules pseudonuls non finis : le noyau de $M \to E$ est donc fini. Ainsi, on a un quasi-isomorphisme

$$a_\Lambda^2(M) \xrightarrow{\sim} a_\Lambda^3(F) \ .$$

Comme la hauteur des idéaux associés à $a_\Lambda^3(F)$ est supérieure à 3, $a_\Lambda^3(F)$ est fini. On en déduit que $a_\Lambda^2(M)$ est fini.

**Définition.** Un $\Lambda$-module de type fini $M$ est dit négligeable si la hauteur des idéaux associés à $M$ est supérieure ou égale à 3. On dit que $M$ vérifie la propriété (A) si $a_\Lambda^i(M)$ est négligeable pour $i \geq 2$.

Un complexe (de longueur fini) de $\Lambda$-modules de type fini est dit suite presque-exacte s'il est exact à des modules négligeables près.

Lorsque $\dim \Lambda = 2$, les modules négligeables sont les modules nuls. Lorsque $\dim \Lambda = 3$, les modules négligeables sont les modules finis. Une suite presque exacte est donc une suite quasi-exacte. Un $\Lambda$-module de type fini $M$ vérifie donc alors la propriété (A) si $a^2(M)$ est fini (cela est automatique pour $a^3(M)$).

Ainsi, si $M$ n'a pas de modules pseudonuls non finis, $M$ vérifie la propriété (A).

Si $\mathfrak{p}$ est un idéal de $\Lambda$, on note $M^{\mathfrak{p}}$ le sous-module des éléments de $M$ annulés par $\mathfrak{p}$.

**1.1.1. Proposition.** *Soit $M$ un $\Lambda$-module de type fini, de torsion et $\mathfrak{p}$ un idéal de $\Lambda$ de hauteur 1.*

*1) Si $M/\mathfrak{p}$ est de $\Lambda/\mathfrak{p}$-torsion, on a la suite exacte naturelle*

$$0 \to a^1_\Lambda(M)/\mathfrak{p} \to a^1_{\Lambda/\mathfrak{p}}(M/\mathfrak{p}) \to a^2_\Lambda(M)^{\mathfrak{p}} \to 0 \ .$$

*2) Si $M$ vérifie la condition (A), on a une suite presque exacte*

$$0 \to a^1_{\Lambda/\mathfrak{p}}(M/\mathfrak{p}) \to a^1_\Lambda(M)/\mathfrak{p} \to a^0_{\Lambda/\mathfrak{p}}(M^{\mathfrak{p}}) \to 0 \ .$$

*Si $M$ est de dimension projective inférieure ou égale à 1, la suite*

$$0 \to a^1_{\Lambda/\mathfrak{p}}(M/\mathfrak{p}) \to a^1_\Lambda(M)/\mathfrak{p} \to a^0_{\Lambda/\mathfrak{p}}(M^{\mathfrak{p}}) \to a^2_{\Lambda/\mathfrak{p}}(M/\mathfrak{p}) \to 0$$

*est exacte.*

*Démonstration.* On déduit de la suite exacte $0 \to \Lambda \xrightarrow{f} \Lambda \to \Lambda/\mathfrak{p} \to 0$ avec $\mathfrak{p} = (f)$ la suite exacte

$$0 \to \operatorname{Ext}^1_\Lambda(M,\Lambda)/\mathfrak{p} \to \operatorname{Ext}^1_\Lambda(M,\Lambda/\mathfrak{p}) \to \operatorname{Ext}^2_\Lambda(M,\Lambda)^{\mathfrak{p}} \to 0 \ .$$

D'autre part, on utilise une résolution de $M$ par des modules de type fini

$$0 \to L' \to L \to M \to 0$$

avec $L$ libre. Si $L''$ est le noyau de $L/\mathfrak{p} \to M/\mathfrak{p}$, on a la suite exacte

$$0 \to M^{\mathfrak{p}} \to L'/\mathfrak{p} \to L'' \to 0 \ .$$

On a le diagramme commutatif dont les lignes et les colonnes sont exactes :

$$\begin{array}{ccccccccc}
 & & & & & & 0 & & \\
 & & & & & & \downarrow & & \\
0 & \to & a^0_{\Lambda/\mathfrak{p}}(M/\mathfrak{p}) & \to & a^0_{\Lambda/\mathfrak{p}}(L/\mathfrak{p}) & \to & a^0_{\Lambda/\mathfrak{p}}(L'') & \to & a^1_{\Lambda/\mathfrak{p}}(M/\mathfrak{p}) \to 0 \\
 & & \| & & \| & & \downarrow & & \downarrow \\
0 & \to & \operatorname{Hom}_\Lambda(M,\Lambda/\mathfrak{p}) & \to & \operatorname{Hom}_\Lambda(L,\Lambda/\mathfrak{p}) & \to & \operatorname{Hom}_\Lambda(L',\Lambda/\mathfrak{p}) & \to & \operatorname{Ext}^1_\Lambda(M,\Lambda/\mathfrak{p}) \to 0 \\
 & & & & & & \downarrow & & \\
 & & & & & & a^0_{\Lambda/\mathfrak{p}}(M^{\mathfrak{p}}) & & \\
 & & & & & & \downarrow & & \\
 & & & & & & a^1_{\Lambda/\mathfrak{p}}(L'') & & \\
 & & & & & & \downarrow & & \\
 & & & & & & a^1_{\Lambda/\mathfrak{p}}(L'/\mathfrak{p}) & &
\end{array}.$$

Comme $a^i_{\Lambda/\mathfrak{p}}(L'') \cong a^{i+1}_{\Lambda/\mathfrak{p}}(M/\mathfrak{p})$ pour $i \geq 1$, on trouve la suite exacte

$$0 \to a^1_{\Lambda/\mathfrak{p}}(M/\mathfrak{p}) \to \operatorname{Ext}^1_\Lambda(M,\Lambda/\mathfrak{p}) \to a^0_{\Lambda/\mathfrak{p}}(M^{\mathfrak{p}}) \to a^2_{\Lambda/\mathfrak{p}}(M/\mathfrak{p}) \ .$$

Lorsque $M/\mathfrak{p}$ est de $\Lambda/\mathfrak{p}$-torsion, il en est de même de $M^{\mathfrak{p}}$, $a^0_{\Lambda/\mathfrak{p}}(M^{\mathfrak{p}})$ est nul et on obtient l'assertion (1). En général, on peut résumer en les deux suites exactes

$$\begin{array}{ccccccccc}
 & & 0 & & & & & & \\
 & & \downarrow & & & & & & \\
 & & a^1_\Lambda(M)/\mathfrak{p} & & & & & & \\
 & & \downarrow & & & & & & \\
0 \to & & a^1_{\Lambda/\mathfrak{p}}(M/\mathfrak{p}) & \to & \operatorname{Ext}^1_\Lambda(M,\Lambda/\mathfrak{p}) & \to & a^0_{\Lambda/\mathfrak{p}}(M^{\mathfrak{p}}) & \to & a^2_{\Lambda/\mathfrak{p}}(M/\mathfrak{p}) \ . \\
 & & & & \downarrow & & & & \\
 & & & & a^2_\Lambda(M)^{\mathfrak{p}} & & & & \\
 & & & & \downarrow & & & & \\
 & & & & 0 & & & &
\end{array}$$

Lorsque $M$ vérifie la propriété (A), $a_\Lambda^2(M)$ est négligeable. Lorsque $M$ est de dimension projective inférieure ou égale à 1, $a_\Lambda^2(M) = 0$, le $\Lambda$-module $L'$ est libre, donc le $\Lambda/\mathfrak{p}$-module $L'/\mathfrak{p}$ est libre et $a_{\Lambda/\mathfrak{p}}^1(L'/\mathfrak{p}) = 0$. D'où la suite exacte de la proposition. □

1.1.2. **Remarque.** (1) L'application $\mathrm{Ext}_\Lambda^1(M, \Lambda/\mathfrak{p}) \to a_{\Lambda/\mathfrak{p}}^0(M^\mathfrak{p})$ dépend du choix d'un générateur $f$ de $\mathfrak{p}$.

(2) Le $\Lambda/\mathfrak{p}$-module $a_\Lambda^1(M)/\mathfrak{p}$ contrôle à la fois la partie libre de $M^\mathfrak{p}$ et la partie de torsion de $M/\mathfrak{p}$. En particulier, si $M/\mathfrak{p}$ est de $\Lambda/\mathfrak{p}$-torsion, $a_\Lambda^1(M)/\mathfrak{p}$ est un module de torsion et on a la suite exacte

$$0 \to a_\Lambda^1(M)/\mathfrak{p} \to a_{\Lambda/\mathfrak{p}}^1(M/\mathfrak{p}) \to a_\Lambda^2(M)^\mathfrak{p} \to 0 \ .$$

Si $M$ est de dimension projective inférieure ou égale à 1, $a_\Lambda^1(M)/\mathfrak{p}$ est égal à $a_{\Lambda/\mathfrak{p}}^1(M/\mathfrak{p})$. Si $M$ vérifie la propriété (A), ($a_\Lambda^2(M)$ est donc fini), le noyau de l'application $a_\Lambda^1(M)/\mathfrak{p} \to a_{\Lambda/\mathfrak{p}}^0(M^\mathfrak{p})$ est le sous-module de torsion de $a_\Lambda^1(M)/\mathfrak{p}$ et est quasi-isomorphe à $a_{\Lambda/\mathfrak{p}}^1(M/\mathfrak{p})$.

1.2. **Construction d'accouplement.** Soit deux $\Lambda$-modules $M$ et $M'$ de $\Lambda$-torsion et un $\Lambda$-homomorphisme $\Theta$

$$M \to a_\Lambda^1(M')$$

autrement dit une application bilinéaire $M \times M' \to \mathrm{Frac}\,\Lambda/\Lambda$. Si $\mathfrak{p}$ est un idéal de $\Lambda$ de hauteur 1, on en déduit un homomorphisme de $\Lambda/\mathfrak{p}$-modules

$$M/\mathfrak{p} \to a_\Lambda^1(M')/\mathfrak{p} \to \mathrm{Ext}_\Lambda^1(M', \Lambda/\mathfrak{p}) \ .$$

L'image du $\Lambda/\mathfrak{p}$-module de torsion $t_{\Lambda/\mathfrak{p}}(M/\mathfrak{p})$ de $M/\mathfrak{p}$ est contenue dans le module de $\Lambda/\mathfrak{p}$-torsion de $\mathrm{Ext}_\Lambda^1(M', \Lambda/\mathfrak{p})$. On a donc le diagramme commutatif dont les lignes sont exactes :

(1.2.1)
$$\begin{array}{ccccccccc}
0 & \to & a_{\Lambda/\mathfrak{p}}^1(M'/\mathfrak{p}) & \to & \mathrm{Ext}_\Lambda^1(M', \Lambda/\mathfrak{p}) & \to & a_{\Lambda/\mathfrak{p}}^0(M'^\mathfrak{p}) & \to & a_{\Lambda/\mathfrak{p}}^2(M'/\mathfrak{p}) \\
 & & & & \uparrow & & & & \\
 & & \uparrow & & a_\Lambda^1(M')/\mathfrak{p} & & & & \\
 & & & & \uparrow & & & & \\
0 & \to & t_{\Lambda/\mathfrak{p}}(M/\mathfrak{p}) & \to & M/\mathfrak{p} & \to & (M/\mathfrak{p})^{**} & \to & \Lambda/\mathfrak{p} - \mathrm{psn}
\end{array}$$

avec $(M/\mathfrak{p})^{**} = \mathrm{Hom}_{\Lambda/\mathfrak{p}}(\mathrm{Hom}_{\Lambda/\mathfrak{p}}(M/\mathfrak{p}, \Lambda/\mathfrak{p}), \Lambda/\mathfrak{p}) = a_{\Lambda/\mathfrak{p}}^0(a_{\Lambda/\mathfrak{p}}^0(M/\mathfrak{p}))$. On obtient ainsi un quasi-homomorphisme de $\Lambda/\mathfrak{p}$-modules

$$\ell_\mathfrak{p}(\Theta) : a_{\Lambda/\mathfrak{p}}^0(a_{\Lambda/\mathfrak{p}}^0(M/\mathfrak{p}) = (M/\mathfrak{p})^{**} \to \mathbb{Q}_p \otimes a_{\Lambda/\mathfrak{p}}^0(M'^\mathfrak{p})$$

et un homomorphisme de $\Lambda/\mathfrak{p}$-modules

$$t_\mathfrak{p}(\Theta) : t_{\Lambda/\mathfrak{p}}(M/\mathfrak{p}) \to a_{\Lambda/\mathfrak{p}}^1(M'/\mathfrak{p}) \to a_{\Lambda/\mathfrak{p}}^1(t_{\Lambda/\mathfrak{p}}(M'/\mathfrak{p})) \ .$$

Contrairement à $t_\mathfrak{p}(\Theta)$, $\ell_\mathfrak{p}(\Theta)$ dépend du choix d'un générateur de $\mathfrak{p}$.

1.2.1. **Lemme.** *Supposons $\Lambda$ de dimension inférieure ou égale à 3. Si $M$ et $M'$ sont des $\Lambda$-modules de torsion vérifiant la propriété (A) et si $\Theta$ est un quasi-isomorphisme de $\Lambda$-modules, $\ell_\mathfrak{p}(\Theta)$ et $t_\mathfrak{p}(\Theta)$ sont des quasi-$\Lambda/\mathfrak{p}$-isomorphismes.*

*Démonstration.* L'hypothèse implique que $\Theta_\mathfrak{p} : M/\mathfrak{p} \to a_\Lambda^1(M')/\mathfrak{p}$ est un quasi-isomorphisme. D'autre part, $a_{\Lambda/\mathfrak{p}}^2(M'/\mathfrak{p})$ est fini ainsi que le conoyau de $M/\mathfrak{p} \to (M/\mathfrak{p})^{**}$. □

Gardons les hypothèses du lemme. On déduit de $t_\mathfrak{p}(\Theta)$ une forme bilinéaire quasi-non dégénérée :

$$t_{\Lambda/\mathfrak{p}}(M/\mathfrak{p}) \times t_{\Lambda/\mathfrak{p}}(M'/\mathfrak{p}) \to \mathrm{Frac}\,\Lambda/\mathfrak{p}/(\Lambda/\mathfrak{p})$$

ou

$$a_{\Lambda/\mathfrak{p}}^1(M/\mathfrak{p}) \times a_{\Lambda/\mathfrak{p}}^1(M'/\mathfrak{p}) \to \mathrm{Frac}\,\Lambda/\mathfrak{p}/(\Lambda/\mathfrak{p}) \ .$$

Lorsqu'on tensorise $\ell_\mathfrak{p}(\Theta)$ par $\mathrm{Frac}\,\Lambda/\mathfrak{p}$, on en déduit un isomorphisme

$$\mathrm{Frac}\,\Lambda/\mathfrak{p} \otimes a_{\Lambda/\mathfrak{p}}^0(a_{\Lambda/\mathfrak{p}}^0(M/\mathfrak{p})) \to \mathrm{Frac}\,\Lambda/\mathfrak{p} \otimes a_{\Lambda/\mathfrak{p}}^0(M'^\mathfrak{p}) \ .$$

En prenant l'inverse et en composant avec l'application induite par $M'^{\mathfrak{p}} \to M'/\mathfrak{p}$, on en déduit une forme bilinéaire

$$a^0_{\Lambda/\mathfrak{p}}(M/\mathfrak{p}) \times a^0_{\Lambda/\mathfrak{p}}(M'/\mathfrak{p}) \to \operatorname{Frac} \Lambda/\mathfrak{p}$$

qui est non dégénérée si et seulement si $\operatorname{Frac} \Lambda/\mathfrak{p} \otimes M'^{\mathfrak{p}} \to \operatorname{Frac} \Lambda/\mathfrak{p} \otimes M'/\mathfrak{p}$ est un isomorphisme, c'est-à-dire si et seulement si le noyau de $M^{\mathfrak{p}} \to M/\mathfrak{p}$ est de $\Lambda/\mathfrak{p}$-torsion.

1.2.1. Lorsque $\Lambda = \mathbb{Z}_p[[\Gamma]]$, $\Lambda$ est muni d'une involution induite par $\tau \to \tau^{-1}$ et que l'on note avec un point. Si $N$ est un $\Lambda$-module, $\dot{N}$ est le module $N$ muni d'une nouvelle action de $\Gamma$ : $\tau \cdot n = \tau^{-1} n$.

Supposons $\Gamma = \mathbb{Z}_p$. Reprenons les suites exactes et la construction : si $\Gamma_n = \Gamma^{p^n}$, les modules $\mathbb{Z}_p[\Gamma/\Gamma_n]$-pseudonuls sont nuls et $a^1_{\mathbb{Z}_p[\Gamma/\Gamma_n]}(\dot{N})$ est égal à $\widehat{t_{\mathbb{Z}_p}(N)}$ muni de l'action usuelle de $\Gamma/\Gamma_n$ sur le dual de Pontryagin : $\tau(f)(n) = f(\tau^{-1}n)$.

Pour $\mathfrak{p}$ l'idéal engendré par $\gamma - 1$, $M^{\mathfrak{p}}$ est le module des invariants $M^{\Gamma}$ et $M/\mathfrak{p}$ le module des coinvariants $M_{\Gamma}$. On obtient un diagramme commutatif dont les lignes sont exactes

$$\begin{array}{ccccccccc}
 & & & & & & \operatorname{Hom}_{\mathbb{Z}_p}(M'_{\Gamma}, \mathbb{Z}_p) & & \\
 & & & & & & \downarrow & & \\
0 & \to & \widehat{t_{\mathbb{Z}_p}(M'_{\Gamma})} & \to & \operatorname{Ext}^1_{\Lambda}(M', \Lambda)_{\Gamma} & \to & \operatorname{Hom}_{\mathbb{Z}_p}(M'^{\Gamma}, \mathbb{Z}_p) & \to & 0 \\
 & & \uparrow & & \uparrow & & \uparrow & & \\
0 & \to & t_{\mathbb{Z}_p}(M_{\Gamma}) & \to & M_{\Gamma} & \to & L_{\mathbb{Z}_p}(M_{\Gamma}) & \to & 0 \ .
\end{array}$$

Autrement dit, on obtient une forme bilinéaire à valeurs dans $\mathbb{Q}_p/\mathbb{Z}_p$

$$\widehat{t_{\mathbb{Z}_p}(M_{\Gamma})} \times \widehat{t_{\mathbb{Z}_p}(M'_{\Gamma})} \to \mathbb{Q}_p/\mathbb{Z}_p$$

et une forme bilinéaire à valeurs dans $\mathbb{Q}_p$

$$M'^*_{\Gamma} \times M^*_{\Gamma} \to \mathbb{Q}_p$$

qui est non dégénérée si $\mathbb{Q}_p \otimes_{\mathbb{Z}_p} M^{\Gamma} \to \mathbb{Q}_p \otimes_{\mathbb{Z}_p} M_{\Gamma}$ est un isomorphisme, En remplaçant $\gamma$ par $\gamma^{p^n}$, on obtient de même une forme sesqui-linéaire

$$\dot{M}'^*_{\Gamma_n} \times M^*_{\Gamma_n} \to \mathbb{Q}_p[\Gamma/\Gamma_n] \ .$$

1.3. **Calcul de l'adjoint.** [19], [16], [1] Prenons $\Lambda = \mathbb{Z}_p[[\Gamma]]$ avec $\Gamma \cong \mathbb{Z}_p^r$. Iwasawa a donné un moyen explicite de calculer $a^1_{\Lambda}(M)$. Soit $\Gamma'$ un sous-groupe isomorphe à $\mathbb{Z}_p$ de $\Gamma$. Soit $\gamma$ un générateur topologique de $\Gamma'$. Posons $\Lambda_n = \mathbb{Z}_p[[\Gamma/\Gamma'_n]]$

1.3.1. **Proposition.** *Soit $M$ un $\Lambda$-module de type fini de torsion tel que $M/(\gamma^{p^n} - 1)$ soit un $\Lambda_n$-module de torsion pour tout entier $n$. Alors*

$$a^1_{\Lambda}(M) = \varprojlim_n a^1_{\Lambda_n}(M/(\gamma^{p^n} - 1)) = \varprojlim_n a^1_{\Lambda_n}(M_{\Gamma'_n})$$

*l'application de transition étant induite par la trace c'est-à-dire par la multiplication par $\sum_{i=0}^{p-1} \gamma^{ip^n}$.*

Cela se déduit des suites exactes

$$0 \to a^1_{\Lambda}(M)_{\Gamma'_n} \to a^1_{\Lambda_n}(M_{\Gamma'_n}) \to a^2_{\Lambda}(M)^{\Gamma'_n}$$

et du fait que la limite projective des $a^2_{\Lambda}(M)^{\Gamma'_n}$ est nulle.

Dans le cas où $\Gamma = \Gamma' = \mathbb{Z}_p$, on obtient le résultat bien connu suivant : si $M_{\Gamma_n}$ est fini pour tout entier $n$,

$$a^1_{\Lambda}(M) \cong \varprojlim_n \widehat{M}^{\Gamma_n} = \varprojlim_n \widehat{M}^{\Gamma_n}/p^n \widehat{M}^{\Gamma_n} \ .$$

Nous renvoyons à [19] ou à [16] pour des précisions et une interprétation en termes de cohomologie locale.

## 2. Modules de Selmer et théorèmes de contrôle

### 2.1. Notations.
Soit $p$ un nombre premier impair, $F$ un corps de nombres, $S$ un nombre fini de places de $F$ contenant les places à l'infini et les places au dessus de $p$. Si $v$ est une place de $F$, on note $G_v$ un groupe de décomposition de $F$ en $v$. Si $L$ est une extension de $F$, on note $S(L)$ l'ensemble des places de $L$ au dessus de $S$. Si $v$ est une place de $F$, la notation $L_v/F_v$ signifie par abus de notation $L_w/F_v$ où $w$ est une place choisie de $L$ au dessus de $v$ (le contexte indiquant que le choix n'a alors pas d'importance).

Soit $V$ une représentation $p$-adique du groupe de Galois absolu $G_F$ de $F$, non ramifiée en dehors de $S$ et ordinaire aux places divisant $p$. Ainsi, si $G_{F,S}$ est le groupe de Galois de la plus grande extension de $F$ non ramifiée en dehors de $S$, $V$ est une représentation $p$-adique de $G_{F,S}$. Soit $T$ un réseau de $V$ stable par $G_F$. On note $V^* = \operatorname{Hom}(V, \mathbb{Q}_p)$, $T^* = \operatorname{Hom}(T, \mathbb{Z}_p)$, $\check{V} = V^*(1)$ le dual de Tate de $V$, $\check{T} = T^*(1)$. Si $W$ est un $\mathbb{Z}_p$-module libre de type fini, on pose $\mathcal{U}(W) = (\mathbb{Q}_p \otimes W)/W$ et $\check{\mathcal{U}}(W) = (\mathbb{Q}_p \otimes \check{W})/\check{W}$.

Nous ferons désormais l'hypothèse suivante :

(Ord) $\qquad\qquad\qquad\qquad V$ est ordinaire aux places divisant $p$.

Dire que $V$ est ordinaire aux places divisant $p$ signifie que pour tout $v \mid p$, il existe une filtration de $G_v$-modules $\operatorname{Fil}_v^j V$ associée à la représentation $p$-adique $V$ telle que le groupe d'inertie en $v$ agit sur le quotient $\operatorname{Fil}_v^j V / \operatorname{Fil}_v^{j+1} V$ par le caractère $\chi^j$ où $\chi$ est le caractère cyclotomique. On pose $\operatorname{Fil}_v^j T = \operatorname{Fil}_v^j V \cap T$.

Soit $\rho$ un caractère continu de $G_F$ à valeurs dans $\mathbb{Z}_p^*$. On note $V \otimes \rho$ la représentation $V$ twistée par le caractère $\rho$. Lorsque $\rho$ est la puissance $k$-ième du caractère cyclotomique, on trouve le twist à la Tate usuel noté $V(k)$. On pose pour simplifier $\mathcal{U}_\rho = \mathcal{U}(T \otimes \rho)$ et $\check{\mathcal{U}}_\rho = \check{\mathcal{U}}(T \otimes \rho) = \check{\mathcal{U}}(T) \otimes \rho^{-1}$.

Nous ferons plus loin l'hypothèse suivante pour certaines extensions $L$ de $F$ :

$(\operatorname{Hyp}(L,V,\rho)) \qquad\qquad (V \otimes \rho)^{G_L} = 0$ et pour $v \mid p$, $((V/\operatorname{Fil}_v^1 V) \otimes \rho)^{G_{L_v}} = 0$.

Enfin, si $F_\infty$ est une $\mathbb{Z}_p$ ou $\mathbb{Z}_p^2$-extension, nous supposerons qu'il n'y a qu'un nombre fini de places au dessus de $p$ dans $F_\infty$.

### 2.2. Groupes de Selmer.
Les modules de Selmer peuvent être définis en prenant la définition de Bloch-Kato ou en prenant celle de Greenberg. Soit

$$H_f^1(F_v, V) = \begin{cases} H^1(G_v/I_v, V^{I_v}) & \text{pour } v \nmid p \\ \operatorname{Im} H^1(F_v, \operatorname{Fil}_v^1 V) \to H^1(F_v, V) \\ = \ker H^1(F_v, V) \to H^1(F_v, V/\operatorname{Fil}_v^1 V) & \text{pour } v \mid p \end{cases}$$

et

$$H_{/f}^1(F_v, V) = H^1(F_v, V)/H_f^1(F_v, V) .$$

On a alors

$$H_{/f}^1(F_v, V) = \begin{cases} H^1(I_v, V)^{G_v/I_v} & \text{pour } v \nmid p \\ H^1(F_v, V/\operatorname{Fil}_v^1 V) & \text{pour } v \mid p \end{cases}$$

(sous l'hypothèse (Hyp(F,V)), en général, l'inclusion seule est vraie pour $v \mid p$). Soit $H_f^1(F_v, T)$ l'image réciproque de $H_f^1(F_v, V)$ dans $H^1(F_v, T)$ et

$$H_{/f}^1(F_v, \mathcal{U}) = H^1(F_v, \mathcal{U})/\mathbb{Q}_p/\mathbb{Z}_p \otimes H_f^1(F_v, T) = H^1(F_v, \mathcal{U})/\operatorname{Im} H_f^1(F_v, V) .$$

On définit $H_f^1(F, T)$ comme le noyau de

$$H^1(G_{F,S}, T) \to \prod_{v \in S} H_{/f}^1(F_v, V) .$$

Ensuite, $H_f^1(F, \mathcal{U}) = H_f^1(F, V/T)$ peut être défini comme le noyau de l'application

$$H^1(G_{F,S}, \mathcal{U}) \to \prod_{v \in S} H_{/f}^1(F_v, \mathcal{U}) .$$

2.3. **Modules d'Iwasawa.** Soit $F_\infty/F$ une $\mathbb{Z}_p$-extension ou une $\mathbb{Z}_p^2$-extension. On pose $\Gamma = \text{Gal}(F_\infty/F)$, $\Lambda = \mathbb{Z}_p[[\text{Gal}(F_\infty/F)]]$ et on note $F_n = F_\infty^{\Gamma^{p^n}}$ le sous-corps de $F_\infty$ fixé par $\Gamma^{p^n}$. Soit alors

$$X_{\infty,f}(F_\infty, \check{T}) = \text{Hom}_{\mathbb{Z}_p}(H^1_f(F_\infty, \mathcal{U}), \mathbb{Q}_p/\mathbb{Z}_p)$$

où $H^1_f(F_\infty, \mathcal{U})$ est la limite inductive des $H^1_f(L, \mathcal{U})$ pour $L$ sous-extension de $F_\infty$. C'est un $\Lambda$-module de type fini. Soit $\rho$ un caractère continu de $G_F$ à valeurs dans $\mathbb{Z}_p^*$ se factorisant par $\Gamma$.

Notons

$$H^1_*(F_v, V \otimes \rho) = \begin{cases} H^1(G_v/I_v, (V \otimes \rho)^{I_v}) & \text{pour } v \nmid p \\ \text{Im } H^1(F_v, \text{Fil}^1_v V \otimes \rho) \to H^1(F_v, V \otimes \rho) \\ = \ker H^1(F_v, V \otimes \rho) \to H^1(F_v, V \otimes \rho/\text{Fil}^1_v V \otimes \rho) & \text{pour } v \mid p \end{cases}$$

et

$$H^1_{/*}(F_v, \mathcal{U}_\rho) = H^1(F_v, \mathcal{U}_\rho)/\text{Im } H^1_*(F_v, V \otimes \rho)$$

Soit $H^1_*(F, \mathcal{U}_\rho)$ le noyau de

$$H^1(G_{S,F}, \mathcal{U}_\rho) \to \prod_{v \in S} H^1_{/*}(F_v, \mathcal{U}_\rho) \subset \prod_{v \mid p} H^1(F_v, \mathcal{U}_\rho)/\text{Im } H^1(F_v, \text{Fil}^1_v V \otimes \rho)$$

$$\prod_{v \in S, v \nmid p} H^1(F_v, \mathcal{U}_\rho)/\text{Im } H^1(G_v/I_v, (V \otimes \rho)^{I_v}) .$$

On définit $H^1_*$, $H^1_{/*}$ comme pour $H^1_f$, $H^1_{/f}$. Lorsque $F_\infty$ contient le corps $L_\rho$ fixé par le noyau de $\rho$, le module

$$X_{\infty,*}(F_\infty, \check{T} \otimes \rho^{-1}) \stackrel{\text{déf}}{=} \text{Hom}_{\mathbb{Z}_p}(H^1_*(F_\infty, \mathcal{U}_\rho), \mathbb{Q}_p/\mathbb{Z}_p) = X_{\infty,f}(F_\infty, \check{T}) \otimes \rho^{-1}$$

est un twist de $X_{\infty,f}(F_\infty, \check{T})$.

2.4. **Théorèmes de contrôle.** On note

$$\mathfrak{Z}(L_v, T \otimes \rho) = \mathfrak{Z}_\rho(L_v) = \begin{cases} \mathcal{U}_\rho(T)^{G_{L_v}}/\text{Im}(V \otimes \rho)^{G_{L_v}} & \text{si } v \nmid p \\ \mathcal{U}_\rho(T/\text{Fil}^1_v T)^{G_{L_v}} & \text{si } v \mid p \end{cases} .$$

Considérons pour $L$ contenu dans $F_\infty$ les applications

$$\Xi_{F_\infty/L}(T \otimes \rho) : H^1(F_\infty/L, \mathcal{U}_\rho^{G_{F_\infty}}) \to \prod_{v \in S(L)} H^1(F_{\infty,v}/L_v, \mathfrak{Z}_\rho(F_{\infty,v})) .$$

Remarquons que seules les places de $F$ qui ne sont pas totalement décomposées dans $F_\infty$ interviennent réellement.

2.4.1. **Proposition** (théorème de contrôle, cas d'une $\mathbb{Z}_p$-extension). *Les applications induites par restriction*

$$H^1_*(F, \mathcal{U}_\rho) \to H^1_*(F_\infty, \mathcal{U}_\rho)^\Gamma = \widehat{X_*(F_\infty, \check{T} \otimes \rho^{-1})}^\Gamma$$

*entrent dans une suite exacte naturelle*

$$0 \to \ker \Xi_{F_\infty/F}(T \otimes \rho) \to H^1_*(F, \mathcal{U}_\rho) \to H^1_*(F_\infty, \mathcal{U}_\rho)^\Gamma$$
$$\to \text{coker } \Xi_{F_\infty/F}(T \otimes \rho) \to \widehat{H^1_*(F, \check{T} \otimes \rho^{-1})} .$$

2.4.2. **Corollaire.** *Sous l'hypothèse* $(\text{Hyp}(F, V, \rho))$, *l'homomorphisme*

$$X_*(F_\infty, \check{T} \otimes \rho^{-1})_\Gamma \to \widehat{H^1_*(F, \mathcal{U}_\rho)}$$

*a un noyau et conoyau finis. Sous l'hypothèse* $(\text{Hyp}(F_\infty, V, \rho))$, *les noyaux et les conoyaux des*

$$X_*(F_\infty, \check{T} \otimes \rho^{-1})_{\Gamma_n} \to \widehat{H^1_*(F_n, \mathcal{U}_\rho)}$$

*sont d'ordre borné par rapport à $n$.*

Il est commode d'introduire un sous-groupe de $H^1(G_{S,F}, \mathcal{U}_\rho)$ un peu plus grand que $H^1_*(F, \mathcal{U}_\rho)$. Il s'agit du noyau de

$$H^1(G_{S,F}, \mathcal{U}_\rho) \to \prod_{v \in S} \tilde{H}^1_{/*}(F_v, \mathcal{U}_\rho)$$

avec

$$\tilde{H}^1_{/*}(F_v, \mathcal{U}_\rho) = \begin{cases} (\prod_{w|v} H^1(F_{\infty,w}, \mathcal{U}_\rho))^\Gamma & \text{si } v \nmid p \\ (\prod_{w|v} H^1(F_{\infty,v}, \mathcal{U}_\rho/\operatorname{Fil}^1_v V \otimes \rho))^\Gamma & \text{si } v \mid p \\ = (\prod_{w|v} H^1(F_{\infty,v}, \mathcal{U}_\rho(T/\operatorname{Fil}^1_v T)))^\Gamma \end{cases}$$

On note aussi $\tilde{H}^1_*(F_v, \mathcal{U}_\rho)$ le noyau de $H^1(F_v, \mathcal{U}_\rho) \to \tilde{H}^1_{/*}(F_v, \mathcal{U}_\rho)$.

L'intérêt d'introduire ce module est le lemme suivant

2.4.3. **Lemme.** *Les suites suivantes sont exactes*

$$0 \to H^1(\Gamma, \mathcal{U}_\rho^{G_{F_\infty}}) \to \tilde{H}^1_*(F, \mathcal{U}_\rho) \to H^1_*(F_\infty, \mathcal{U}_\rho)^\Gamma \to 0$$

$$0 \to X_*(F_\infty, \check{T} \otimes \rho^{-1})_\Gamma \to \widehat{\tilde{H}^1_*(F, \mathcal{U}_\rho)} \to \widehat{H^1(\Gamma, \mathcal{U}_\rho^{G_{F_\infty}})} \to 0$$

*Démonstration.* On a en effet le diagramme commutatif dont les lignes et les colonnes sont exactes

$$\begin{array}{ccccccc}
& & & & 0 & & \\
& & & & \uparrow & & \\
0 & \to & H^1_*(F_\infty, \mathcal{U}_\rho)^\Gamma & \to & H^1(G_{S,F_\infty}, \mathcal{U}_\rho)^\Gamma & \to & \left(\prod_{v \in S(F_\infty)} H^1_{/*}(F_{\infty,v}, \mathcal{U}_\rho)\right)^\Gamma \\
& & \uparrow & & \uparrow & & \uparrow \\
0 & \to & \tilde{H}^1_*(F, \mathcal{U}_\rho) & \to & H^1(G_{S,F}, \mathcal{U}_\rho) & \to & \prod_{v \in S} \tilde{H}^1_{/*}(F_v, \mathcal{U}_\rho) \\
& & & & \uparrow & & \\
& & & & H^1(\Gamma, \mathcal{U}_\rho^{G_{F_\infty}}) & & \\
& & & & \uparrow & & \\
& & & & 0 & &
\end{array}$$

Il s'agit de voir que la flèche verticale de droite est injective. Lorsque $v \nmid p$ est totalement décomposée dans $F_\infty$, $H^1_{/*}(F_v, \mathcal{U}_\rho) = (\prod_{w|v} H^1_{/*}(F_{\infty,v}, \mathcal{U}_\rho))^\Gamma$ car $\Gamma$ agit simplement par permutation des facteurs. Soit $v$ ne divisant pas $p$ et non totalement décomposée dans $F_\infty$. Si $w$ est une place de $F_\infty$ au dessus de $v$, l'extension $F_{\infty,w}$ est l'unique extension de $F_v$ non ramifiée et de groupe de Galois un pro-$p$-groupe. Donc, $H^1_{/*}(F_{\infty,w}, \mathcal{U}_\rho)$ est égal à $H^1(F_{\infty,w}, \mathcal{U}_\rho)^\Gamma$. On en déduit que $\tilde{H}^1_{/*}(F_v, \mathcal{U}_\rho) \to \left(\prod_{w|v} H^1_{/*}(F_{\infty,v}, \mathcal{U}_\rho)\right)^\Gamma$ est un isomorphisme. Lorsque $v \mid p$, l'assertion est claire. □

*Démonstration de la proposition.* La différence entre $\tilde{H}^1_*(F, \mathcal{U}_\rho)$ et $H^1_*(F, \mathcal{U}_\rho)$ est calculée par la suite exacte suivante, conséquence de la suite exacte de Poitou-Tate :

$$0 \to H^1_*(F, \mathcal{U}_\rho) \to \tilde{H}^1_*(F, \mathcal{U}_\rho) \to \prod_{v \in S} \tilde{H}^1_*(F_v, \mathcal{U}_\rho)/H^1_*(F_v, \mathcal{U}_\rho) \to \widehat{H^1_*(F, \check{T} \otimes \rho^{-1})}$$

et il n'est pas difficile de voir que $\prod_{v \in S} \tilde{H}^1_*(F_v, \mathcal{U}_\rho)/H^1_*(F_v, \mathcal{U}_\rho)$ est exactement l'ensemble d'arrivée de $\Xi_{F_\infty/F}(T \otimes \rho)$. En effet, cela est clair pour la contribution des places totalement décomposées dans $F_\infty$. Pour une place $v$ non totalement décomposée dans $F_\infty$ et ne divisant pas $p$, on a d'après le calcul précédent

$$\tilde{H}^1_*(F_v, \mathcal{U}_\rho) = H^1(F_{\infty,v}/F_v, \mathcal{U}_\rho^{G_{F_\infty,v}})$$
$$H^1_*(F_v, \mathcal{U}_\rho) = \mathbb{Q}_p/\mathbb{Z}_p \otimes H^1(F_{\infty,v}/F_v, (T \otimes \rho)^{G_{F_\infty,v}})$$
$$= H^1(F_{\infty,v}/F_v, \mathcal{U}_\rho(T^{G_{F_\infty,v}}))$$

car $\operatorname{Gal}(F_{\infty,v}/F_v)$ est de dimension cohomologique 1. Donc,

$$\tilde{H}^1_*(F_v, \mathcal{U}_\rho)/H^1_*(F_v, \mathcal{U}_\rho) = H^1(F_{\infty,v}/F_v, \mathfrak{Z}_\rho(F_{\infty,v})) \ .$$

Remarquons que par la dualité de Tate, c'est aussi le dual de Pontryagin du sous-$\mathbb{Z}_p$-module de torsion de $H^1(F_{\infty,v}, \check{T} \otimes \rho^{-1})^{G_{F_\infty,v}}$ dont le cardinal est le nombre de Tamagawa local en $v$ de $\check{T} \otimes \rho^{-1}$.

Soit maintenant une place $v$ divisant $p$. On a par définition le diagramme commutatif et exact suivant

$$\begin{array}{ccccccccc}
& & & & 0 & & 0 & & \\
& & & & \uparrow & & \uparrow & & \\
0 & \to & \tilde{H}^1_*(F_v,\mathcal{U}_\rho) & \to & H^1(F_v,\mathcal{U}_\rho) & \to & H^1(F_{\infty,v},\mathcal{U}_\rho(T/\operatorname{Fil}^1_v T))^{\Gamma_v} & & \\
& & \uparrow & & \uparrow & & \uparrow & & \\
0 & \to & H^1_*(F_v,\mathcal{U}_\rho) & \to & H^1(F_v,\mathcal{U}_\rho) & \to & H^1(F_v,\mathcal{U}_\rho(T/\operatorname{Fil}^1_v T)) & \to & H^2 \\
& & \uparrow & & \uparrow & & \uparrow & & \\
& & 0 & & 0 & & H^1(F_{\infty,v}/F_v, (\mathcal{U}_\rho(T/\operatorname{Fil}^1_v T))^{G_{F_{\infty,v}}}) & & \\
& & & & & & \uparrow & & \\
& & & & & & 0 & &
\end{array}$$

L'image de $H^1(F_{\infty,v}/F_v,\mathcal{U}_\rho(T/\operatorname{Fil}^1_v T)^{G_{F_{\infty,v}}})$ dans $H^2 = H^2(F_v,\mathcal{U}_\rho(\operatorname{Fil}^1_v T))$ est nulle. D'où l'assertion sur la contribution en $p$. □

Pour démontrer le corollaire, on remarque que si $v \nmid p$ et que $v$ n'est pas totalement décomposée dans $F_\infty$, $H^1(F_{\infty,v}/F_v,\mathcal{U}_\rho^{G_{F_{\infty,v}}}/(V \otimes \rho)^{G_{F_{\infty,v}}})$ est dual du sous-groupe de torsion de $H^1(I_v,\check{T}{\otimes}\rho^{-1})^{G_v/I_v}$ et a comme cardinal le nombre de Tamagawa local $\operatorname{Tam}_v(\check{T}{\otimes}\rho^{-1})$. Ainsi, il vaut 0 si $V$ a bonne réduction en $v$. Lorsqu'on remplace $F$ par $F_n$, $\operatorname{Tam}_{F_n,v}(\check{T}{\otimes}\rho^{-1})$ est borné par rapport à $n$ ([26, 2.2.6]). Pour $v \mid p$, l'hypothèse (Hyp($F_\infty, V, \rho$)) implique que $H^1(F_{\infty,v}/F_{n,v},\mathcal{U}_\rho^{G_{F_{\infty,v}}})$ est fini et d'ordre borné par rapport à $n$.

### 2.5. Théorème de contrôle : cas d'une $\mathbb{Z}_p^2$-extension.
On suppose maintenant que $F_\infty$ est une $\mathbb{Z}_p^2$-extension de $F$. On a les théorèmes de contrôle suivants relatifs à la descente de $F_\infty$ à une sous-$\mathbb{Z}_p$-extension $L_\infty$ (on note alors $\Lambda_{L_\infty} = \mathbb{Z}_p[[\operatorname{Gal}(L_\infty/F)]]$). Notons $\Xi_{F_\infty/L_\infty}(T \otimes \rho)$ l'application

$$H^1(F_\infty/L_\infty, \mathcal{U}_\rho^{G_{F_\infty}}) \to \prod_{v \in S(L_\infty)} H^1(F_{\infty,v}/L_{\infty,v}, \mathfrak{Z}_\rho(F_{\infty,v})) \ .$$

Seules les places totalement décomposées dans $L_\infty$ et les places divisant $p$ interviennent en fait. En effet, dans le cas contraire, elles sont nécessairement totalement décomposées dans $F_\infty/L_\infty$. Notons enfin $\check{H}^1_*(L_\infty,\check{T}{\otimes}\rho^{-1})$ la limite projective des $H^1_*(L_n,\check{T}{\otimes}\rho^{-1})$ pour les applications de corestriction.

**2.5.1. Proposition** (théorème de contrôle, cas d'une $\mathbb{Z}_p^2$-extension). *Soit $L_\infty$ une sous-$\mathbb{Z}_p$-extension de $F_\infty$. L'application de restriction induit par dualité un homomorphisme*

$$r_{F_\infty/L_\infty} : X_{\infty,*}(F_\infty,\check{T}{\otimes}\rho^{-1})_{\operatorname{Gal}(F_\infty/L_\infty)} \to X_{\infty,*}(L_\infty,\check{T}{\otimes}\rho^{-1})$$

*qui se trouve dans une suite exacte naturelle de $\Lambda_{L_\infty}$-modules*

$$\check{H}^1_*(L_\infty,\check{T}{\otimes}\rho^{-1}) \to \operatorname{coker}\widehat{\Xi_{F_\infty/L_\infty}}(T \otimes \rho) \to X_{\infty,*}(F_\infty,\check{T}{\otimes}\rho^{-1})_{\operatorname{Gal}(F_\infty/L_\infty)}$$

$$\to X_{\infty,*}(L_\infty,\check{T}{\otimes}\rho^{-1}) \to \ker\widehat{\Xi_{F_\infty/L_\infty}}(T \otimes \rho) \to 0 \ .$$

*Sous l'hypothèse (Hyp($L_\infty,V,\rho$)), le noyau de $r_{F_\infty/L_\infty}$ est fini et son conoyau est annulé par une puissance de $p$. Lorsqu'il y a un nombre fini de places au dessus de $S$ dans $L_\infty$ et sous (Hyp($F_\infty,V,\rho$)), les noyaux et conoyaux des $r_{F_\infty/F_nL_\infty}$ sont finis et d'ordre borné par rapport à $n$.*

*Démonstration.* La suite exacte se démontre en utilisant le diagramme exact et commutatif suivant avec $\Gamma' = \operatorname{Gal}(F_\infty/L_\infty)$,

$$\begin{array}{ccccccc}
& & & & 0 & & 0 \\
& & & & \uparrow & & \uparrow \\
0 \to & H^1_*(F_\infty,\mathcal{U}_\rho)^{\Gamma'} & \to & H^1(G_{F_\infty,S},\mathcal{U}_\rho)^{\Gamma'} & \to & (\prod_{w \in S(F_\infty)} H^1_{/*}(F_{\infty,w},\mathcal{U}_\rho))^{\Gamma'} & \\
& \uparrow & & \uparrow & & \uparrow & \\
0 \to & H^1_*(L_\infty,\mathcal{U}_\rho) & \to & H^1(G_{L_\infty,S},\mathcal{U}_\rho) & \to & \prod_{v \in S(L_\infty)} H^1_{/*}(L_{\infty,v},\mathcal{U}_\rho) & \\
& & & \uparrow & & \uparrow & \\
& & & H^1(\Gamma',\mathcal{U}_\rho^{G_{F_\infty}}) & \to & \prod_{v \in S(L_\infty)} H^1(F_{\infty,v}/L_{\infty,v},\mathfrak{Z}_\rho(F_{\infty,v})) & \\
& & & \uparrow & & \uparrow & \\
& & & 0 & & 0 &
\end{array}$$

En effet, soit $w$ une place de $F_\infty$ ne divisant pas $p$ et $v$ sa restriction à $L_\infty$. Si $w$ n'est pas totalement décomposée dans $F_\infty$, on a $H^1_{/*}(F_{\infty,w}, \mathcal{U}_\rho) = H^1(F_{\infty,w}, \mathcal{U}_\rho)$ car $H^1_*(F_{\infty,w}, \mathcal{U}_\rho) = 0$. Si $v$ n'est pas totalement décomposée non plus dans $L_\infty$, on a alors aussi $H^1_{/*}(L_{\infty,w}, \mathcal{U}_\rho) = H^1(L_{\infty,w}, \mathcal{U}_\rho)$ et l'application d'inflation est un isomorphisme. Si $v$ est totalement décomposée dans $L_\infty$, le noyau de $H^1_{/*}(L_{\infty,w}, \mathcal{U}_\rho) \to \prod_{w|v} H^1_{/*}(F_{\infty,w}, \mathcal{U}_\rho)$ est égal au quotient

$$H^1(F_{\infty,v}/L_{\infty,v}, \mathcal{U}_\rho)/\mathbb{Q}_p/\mathbb{Z}_p \otimes H^1(F_{\infty,v}/L_{\infty,v}, (T \otimes \rho)^{G_{F_{\infty,v}}}) \cong H^1(F_{\infty,v}/L_{\infty,v}, \mathfrak{Z}_\rho(F_{\infty,v})) \ .$$

Si $w$ est totalement décomposée dans $F_\infty$, l'assertion est triviale. Si $w$ est une place divisant $p$, le noyau de l'application inflation est isomorphe à $H^1(F_{\infty,v}/L_{\infty,v}, \mathfrak{Z}_\rho(F_{\infty,v}))$. Cela démontre les assertions sur le diagramme précédent. Par une des variantes de la suite exacte de Poitou-Tate, le conoyau de $H^1(G_{L_{\infty,S}}, \mathcal{U}_\rho) \to \prod_{v \in S(L_\infty)} H^1_{/*}(L_{\infty,v}, \mathcal{U}_\rho)$ est contenu dans le dual de Pontryagin de $\check{H}^1_*(L_\infty, \check{T} \otimes \rho^{-1})$. On en déduit la proposition. □

Un cas particulier est le cas où $\rho$ est le caractère trivial.

2.5.2. **Corollaire.** *Soit $L_\infty$ une sous-$\mathbb{Z}_p$-extension de $F_\infty$. Il existe une suite exacte naturelle de $\Lambda_{L_\infty}$-modules*

$$\check{H}^1_f(L_\infty, \check{T}) \to \operatorname{coker} \widehat{\Xi_{F_\infty/L_\infty}}(T)$$
$$\to X_{\infty, f}(F_\infty, \check{T})_{\operatorname{Gal}(F_\infty/L_\infty)} \to X_{\infty, f}(L_\infty, \check{T}) \to \ker \widehat{\Xi_{F_\infty/L_\infty}}(T) \to 0 \ .$$

2.5.3. **Corollaire.** *Soit $L_\infty$ une sous-$\mathbb{Z}_p$-extension de $F_\infty$. Si $X_{\infty, *}(L_\infty, \check{T} \otimes \rho^{-1})$ est un $\Lambda_{L_\infty}$-module de torsion, alors $X_{\infty, *}(F_\infty, \check{T} \otimes \rho^{-1})$ est un $\Lambda$-module de torsion.*

2.6. **Construction de l'adjoint.**

2.6.1. Supposons que $F_\infty$ est une $\mathbb{Z}_p$-extension de $F$. Nous dirons que $\rho$ est admissible (pour $V$ et $F_\infty$) si pour tout entier $n$, les $X_*(F_\infty, \check{T} \otimes \rho^{-1})_{\Gamma_n}$ sont finis. Si $\rho$ est admissible, nécessairement $X_{\infty, *}(F_\infty, \check{T} \otimes \rho^{-1})$ est de $\Lambda$-torsion.

Il existe un caractère de $\operatorname{Gal}(F_\infty/F)$ dans $\mathbb{Z}_p^*$ admissible pour $V$ si et seulement si

$(\operatorname{Tors}(F_\infty, V))$ $\qquad X_f(F_\infty, \check{T})$ est un $\Lambda$-module de torsion

et il en existe alors une infinité. En effet, fixons un caractère non trivial de $\operatorname{Gal}(F_\infty/F)$ dans $\mathbb{Z}_p^*$ ; si $M$ est un $\Lambda$-module de type fini et de torsion, pour tout entier $k$ sauf un nombre fini, $(M \otimes \rho^k)_{\Gamma_n}$ est de torsion pour tout entier $n$. En effet, si $H$ est une série caractéristique de $M$ (en particulier $H$ annule $M$), $(M \otimes \rho^k)_{\Gamma_n}$ est fini si et seulement si $H(u^k \zeta_n - 1)$ est non nul pour $\zeta_n$ une racine de l'unité d'ordre $p^n$ et $u = \rho(\gamma)$. Comme $H$ n'a qu'un nombre fini de zéros par le théorème de préparation de Weierstrass, le fait précédent s'en déduit.

Les applications

$$H^1_*(F_n, \mathcal{U}_\rho) \to X_*(\widehat{F_\infty, \check{T}} \otimes \rho^{-1})^{\Gamma_n}$$

induisent par passage à la limite projective pour les applications de corestriction un $\Lambda$-homomorphisme

$$\mathcal{A}^{(\rho)}_{F_\infty} : \varprojlim_n H^1_*(F_n, \mathcal{U}_\rho) \to a^1_\Lambda(\dot{X}_{\infty, *}(F_\infty, \check{T} \otimes \rho^{-1}))$$

(proposition 1.3.1). Soit

$$\xi_{F_\infty}(T \otimes \rho) : \mathcal{U}_\rho^{G_{F_\infty}} \to \prod_{v \in S^{nd}(F_\infty)} \mathfrak{Z}_\rho(F_{\infty,v})$$

où $S^{nd}(F_\infty)$ désigne l'ensemble des places de $S(F_\infty)$ non totalement décomposées sur $F$.

2.6.1. **Proposition.** *Soit $F_\infty/F$ une $\mathbb{Z}_p$-extension telle que $(\operatorname{Tors}(F_\infty, V))$ soit vérifiée et soit $\rho$ admissible pour $F_\infty$. On a la suite exacte naturelle*

$$0 \to \ker \xi_{F_\infty}(T \otimes \rho) \to \varprojlim_n H^1_*(F_n, \mathcal{U}_\rho) \stackrel{\mathcal{A}^{(\rho)}_{F_\infty}}{\to} a^1_\Lambda(\dot{X}_{\infty, *}(F_\infty, \check{T} \otimes \rho^{-1}))$$
$$\to \operatorname{coker} \xi_{F_\infty}(T \otimes \rho) \to H^1_*(\widehat{F_\infty, \check{T}} \otimes \rho^{-1}) \ .$$

**2.6.2. Corollaire.** *Sous les hypothèses de la proposition, si de plus* $(\mathrm{Hyp}(F_\infty, V, \rho))$ *est vérifié,* $\mathcal{A}^{(\rho)}_{F_\infty}$ *est un quasi-isomorphisme.*

Le corollaire se déduit de la finitude du noyau et du conoyau de $\xi_{F_\infty}(T \otimes \rho)$ (il n'y a qu'un nombre fini de places dans $S^{nd}(F_\infty)$).

*Démonstration de la proposition.* La proposition se déduit de la proposition 1.3.1 et de la proposition 2.4.1 (remarquons que l'application de corestriction devient dans l'isomorphisme $H^1(\Gamma, S) \cong S_\Gamma$ l'application induite par l'identité sur $S$). □

**2.6.3. Remarque.** Supposons de plus que $\rho^{-1}$ est admissible pour $\check{T}$ et pour $F_\infty$. Alors, $H^1_*(F_n, \check{T} \otimes \rho^{-1})$ est fini pour tout entier $n$ et est égal au sous-groupe de $\mathbb{Z}_p$-torsion de $H^1_*(F_n, \check{T} \otimes \rho^{-1})$, c'est-à-dire à $\check{\mathcal{U}}_\rho^{G_{F_n}}$. Donc, sous cette hypothèse,

$$H^1_*(F_\infty, \check{T} \otimes \rho^{-1}) = \check{\mathcal{U}}_\rho^{G_{F_\infty}} = (\check{V} \otimes \rho^{-1}/\check{T} \otimes \rho^{-1})^{G_{F_\infty}}.$$

**2.6.2.** Prenons maintenant pour $F_\infty/F$ une $\mathbb{Z}_p^2$-extension. Soit $F_n$ le corps fixe par $\Gamma_n = \Gamma^{p^n}$.

Si $L_\infty$ est une sous-$\mathbb{Z}_p$-extension de $F_\infty/F$, on note $L_{\infty,n} = L_\infty F_n$. Ainsi, $L_{\infty,n+1}/L_{\infty,n}$ est une extension d'ordre $p$ (pour $n$ assez grand). D'autre part, fixons une $\mathbb{Z}_p$-extension $L'_\infty$ de $F$ telle que $F_\infty = L_\infty L'_\infty$. On pose $\Lambda_{L_{\infty,n}} = \mathbb{Z}_p[[\mathrm{Gal}(L_{\infty,n}/L'_n)]]$. Si $M$ est un $\mathbb{Z}_p[[\mathrm{Gal}(L_{\infty,n}/F)]]$-module, $a^1_{\mathbb{Z}_p[[\mathrm{Gal}(L_{\infty,n}/F)]]}(M)$ et $a^1_{\Lambda_{L_{\infty,n}}}(M)$ muni de sa structure naturelle de $\mathbb{Z}_p[[\mathrm{Gal}(L_{\infty,n}/F)]]$-modules s'identifient canoniquement ([16]). La norme de $L_{\infty,n+1}$ à $L_{\infty,n}$ induit alors des homomorphismes naturels :

$$a^1_{\mathbb{Z}_p[[\mathrm{Gal}(L_{\infty,n+1}/F)]]}(M_{\Gamma_{n+1}}) \to a^1_{\mathbb{Z}_p[[\mathrm{Gal}(L_{\infty,n}/F)]]}(M_{\Gamma_n}).$$

Choisissons une $\mathbb{Z}_p$-extension $L_\infty$ telle que $X_{\infty,*}(F_\infty, \check{T} \otimes \rho^{-1})_{\mathrm{Gal}(F_\infty/L_{\infty,n})}$ soit de $\Lambda_{L_{\infty,n}}$-torsion pour tout entier $n$. Par la proposition 2.5.1, cela est équivalent à ce que $X_{\infty,*}(L_{\infty,n}, \check{T} \otimes \rho^{-1})$ soit de $\Lambda_{L_{\infty,n}}$-torsion pour tout entier $n$. On dit alors que $\rho$ est admissible pour $F_\infty/L_\infty$.

**2.6.4. Proposition.** *On suppose vérifiés* $(\mathrm{Hyp}(F_\infty, V, \rho))$, *que*

$(\mathrm{Tors}(F_\infty, V, \rho))$     $X_{\infty,*}(F_\infty, \check{T} \otimes \rho^{-1})$ *est un $\Lambda$-module de torsion*

*et que $\rho$ est admissible pour $F_\infty/L_\infty$. Les applications naturelles*

$$r_n : a^1_{\Lambda_{L_{\infty,n}}}(X_{\infty,*}(L_{\infty,n}, \check{T} \otimes \rho^{-1})) \to a^1_{\Lambda_{L_{\infty,n}}}(X_{\infty,*}(F_\infty, \check{T} \otimes \rho^{-1})_{\mathrm{Gal}(F_\infty/L_{\infty,n})})$$

*induisent un $\Lambda$-homomorphisme $r_\infty$ injectif*

$$\varprojlim_n a^1_{\Lambda_{L_{\infty,n}}}(X_{\infty,*}(L_{\infty,n}, \check{T} \otimes \rho^{-1})) \to a^1_\Lambda(X_{\infty,*}(F_\infty, \check{T} \otimes \rho^{-1}))$$

*et on a la suite quasi-exacte*

$$0 \to \varprojlim_n a^1_{\Lambda_{L_{\infty,n}}}(X_{\infty,*}(L_{\infty,n}, \check{T} \otimes \rho^{-1})) \to a^1_\Lambda(X_{\infty,*}(F_\infty, \check{T} \otimes \rho^{-1}))$$

$$\to \prod_{v \in S^{nd}(F_\infty/L_\infty)} \mathfrak{Z}_\rho(F_{\infty,v})$$

*Démonstration.* Posons $\Lambda_n = \Lambda_{L_{\infty,n}}$, $\Gamma'_n = \mathrm{Gal}(F_\infty/L_{\infty,n})$ et $a^1_n = a^1_{\Lambda_n}$, $M_n = X_{\infty,*}(L_{\infty,n}, \check{T} \otimes \rho^{-1})$, $M = X_{\infty,*}(F_\infty, \check{T} \otimes \rho^{-1})$, $\mathfrak{Z} = \prod_{v \in S(F_\infty)} \mathfrak{Z}_\rho(F_{\infty,v})$, $\Xi_n = \Xi_{F_\infty/L_{\infty,n}}(T \otimes \rho)$. La suite exacte de la proposition 2.5.1 appliquée à $F_\infty/L_{\infty,n}$ devient

$$\widehat{\mathrm{coker}\,\Xi_n} \to M_{\Gamma'_n} \to M_n \to \widehat{\ker \Xi_n} \to 0$$

et on a la suite exacte tautologique

$$0 \to \widehat{\mathrm{coker}\,\Xi_n} \to \widehat{H^1(\Gamma'_n, \mathfrak{Z})} \to \widehat{H^1(\Gamma'_n, \mathcal{U}_\rho^{G_{F_\infty}})} \to \widehat{\ker \Xi_n} \to 0.$$

Soit $M'_n$ l'image de $M_{\Gamma'_n}$ dans $M_n$. On a alors les suites exactes

$$0 \to a^1_n(M'_n) \to a^1_n(M_{\Gamma'_n}) \to a^1_n(\widehat{\mathrm{coker}\,\Xi_n})$$
$$0 \to a^1_n(M_n) \to a^1_n(M'_n) \to a^2_n(\widehat{\ker \Xi_n}).$$

On en déduit l'injectivité de $a_n^1(M_n) \to a_n^1(M_{\Gamma'_n})$ et par passage à la limite celle de $\varprojlim_n a_n^1(M_n) \to a_\Lambda^1(M_{\Gamma'_n})$.

D'autre part, on déduit de la suite exacte tautologique la suite exacte

$$0 \to a_n^1(\widehat{H^1(\Gamma'_n, \mathfrak{Z})}) \to a_n^1(\widehat{\operatorname{coker} \Xi_n}) \to R_n \to 0$$

avec $R_n$ d'ordre borné par rapport à $n$ (on utilise le fait que $a_n^1(R) = 0$ si $R$ est un module fini).

Nous allons maintenant raisonner à des modules finis près d'ordre borné par rapport à $n$ (on parle alors de suites quasi-exactes et de quasi-isomorphismes contrôlés) : on a la suite quasi-exacte contrôlée :

$$0 \to a_n^1(M_n) \to a_n^1(M_{\Gamma'_n}) \to a_n^1(\widehat{\operatorname{coker} \Xi_n})$$

et le quasi-isomorphisme contrôlé

$$a_n^1(\widehat{H^1(\Gamma'_n, Z)}) \cong a_n^1(\widehat{\operatorname{coker} \Xi_n})$$

Comme $\mathfrak{Z}$ est annulé par une puissance de $p$, $a_n^1(\widehat{H^1(\Gamma'_n, \mathfrak{Z})}) \cong H^1(\Gamma'_n, \mathfrak{Z})$ et la limite projective des $a_n^1(\widehat{\operatorname{coker} \Xi_n})$ est quasi-isomorphe à $\mathfrak{Z}$. La proposition se déduit alors de la proposition 1.3.1 □

**2.6.5. Remarque.** On peut être plus précis sous une hypothèse dont on montrera plus tard qu'elle est vraie. Supposons que le plus grand sous-$\Lambda_n$-module fini de $M_n = X_{\infty,*}(L_{\infty,n}, \check{T} \otimes \rho^{-1})$ est d'ordre borné par rapport à $n$. Alors la dernière flèche est quasi-surjective. En effet, comme $M'_n$ est contenu dans $M_n$, $a_n^2(M'_n)$ est fini d'ordre borné par rapport à $n$. On a donc la suite quasi-exacte

$$0 \to \varprojlim_n a_{\Lambda_{L_\infty,n}}^1(X_{\infty,*}(L_{\infty,n}, \check{T} \otimes \rho^{-1})) \to a_\Lambda^1(X_{\infty,*}(F_\infty, \check{T} \otimes \rho^{-1}))$$
(2.6.1)
$$\to \prod_{v \in S^{nd}(F_\infty/L_\infty)} \mathfrak{Z}_\rho(F_{\infty,v}) \to 0$$

3. Construction d'accouplements entre modules de Selmer

3.1. **Accouplements de Cassels-Tate.**

3.1.1. **Théorème** (Flach). *Il existe un homomorphisme naturel*

$$\operatorname{Cassels}_F(T \otimes \rho) : H_*^1(F, \check{\mathcal{U}}_\rho) \times H_*^1(F, \mathcal{U}_\rho) \to \mathbb{Q}_p/\mathbb{Z}_p$$

*qui induit un isomorphisme*

$$C_F(T \otimes \rho) : \widehat{H_*^1(F, \check{\mathcal{U}}_\rho)} / \operatorname{div} \to H_*^1(F, \mathcal{U}_\rho) / \operatorname{div}$$

*où $M/\operatorname{div}$ désigne le quotient d'un $\mathbb{Z}_p$-module $M$ par sa partie divisible. En particulier, si $H_*^1(F, \check{\mathcal{U}}_\rho)$ et $H_*^1(F, \mathcal{U}_\rho)$ sont finis, on en déduit un isomorphisme*

$$C_F(T \otimes \rho) : \widehat{H_*^1(F, \check{\mathcal{U}}_\rho)} \to H_*^1(F, \mathcal{U}_\rho) \ .$$

*On a les propriétés suivantes :*

(1) *Si $L/F$ est une extension finie,*

$$\operatorname{Cassels}_F(T \otimes \rho)(x, \operatorname{cores}_{L/F} y) = \operatorname{Cassels}_L(\operatorname{res}_{L/F} x, y)$$

(2) *Si $F/F_1$ est une extension galoisienne et $\sigma \in \operatorname{Gal}(F/F_1)$,*

$$\operatorname{Cassels}_F(T \otimes \rho)(\sigma x, \sigma y) = \operatorname{Cassels}_F(T \otimes \rho)(x, y) \ .$$

(3) *Soit $L$ un corps contenant le corps fixe par le noyau de $\rho^{p^n}$. Soit $x \in H_*^1(L, \check{\mathcal{U}}_{p^n})$ et $y \in H_*^1(L, \mathcal{U}_{p^n})$. Alors, si $Tw_\rho(x)$ (resp. $Tw_{\rho^{-1}}(x)$ désigne le $\rho$-ième twist de $x$ (resp. le $\rho^{-1}$-ième twist de $y$), on a*

$$\operatorname{Cassels}_L(T \otimes \rho)(Tw_\rho(x), Tw_{\rho^{-1}}(x)) = \operatorname{Cassels}_L(T)(x, y)$$

(4) *Le dual de $C_F(T \otimes \rho)$ par la dualité de Pontryagin est $C_F(\check{T} \otimes \rho^{-1})$ .*

**3.1.2. Remarque.** 1) La partie divisible de $H^1_*(F, \mathcal{U}_\rho)$ est $\mathbb{Q}_p/\mathbb{Z}_p \otimes H^1_*(F, T \otimes \rho)$.

2) Pour $\rho$ le caractère trivial, $V = V_p(E)$ et en utilisant l'accouplement de Weil pour identifier $V \to \check{V} = V^*(1)$, l'accouplement obtenu est l'accouplement de Cassels. L'accouplement de Weil étant alterné, l'accouplement de Cassels est une forme bilinéaire alternée. C'est ce qui permet de montrer que l'ordre du quotient du groupe de Tate-Shafarevich par sa partie divisible est un carré.

La démonstration du théorème 3.1.1 est faite dans [9], les deux sous-espaces de $V \otimes \rho$ et $\check{V} \otimes \rho^{-1}$ que sont $\mathrm{Fil}^1_v V \otimes \rho$ et $\mathrm{Fil}^1_v \check{V} \otimes \rho^{-1}$ sont orthogonaux dans la dualité naturelle $V \otimes \rho \times \check{V} \otimes \rho^{-1} \to \mathbb{Q}_p(1)$ (voir aussi [11, §5.4]) Pour le comportement par twist, il suffit de reprendre la définition en remarquant que pour $\tau \in G_L$, $\rho(\tau) \equiv 1 \mod p^n$. Les différentes cochaines construites diffèrent alors d'éléments de $T$ et finalement l'image est la même dans $\frac{1}{p^n}\mathbb{Z}/\mathbb{Z}$.

## 3.2. Dualité : cas d'une $\mathbb{Z}_p$-extension.

**3.2.1. Théorème.** *Soit $F_\infty/F$ une $\mathbb{Z}_p$-extension. Soit $\rho$ un caractère continu de $G_F$ à valeurs dans $\mathbb{Z}_p^*$ tel que $(\mathrm{Tors}(F_\infty, V, \rho))$ soit vérifiée et tel que $\rho$ soit admissible pour $F_\infty$ et $V$. Les applications $C_{F_n}(T \otimes \rho)$ induisent un $\Lambda$-homomorphisme quasi-injectif*

$$C_{F_\infty}(T \otimes \rho) : X_{\infty,*}(F_\infty, T \otimes \rho) \to a^1_\Lambda(\dot{X}_{\infty,*}(F_\infty, \check{T} \otimes \rho^{-1}))$$

*et on a plus précisément la suite exacte*

$$0 \to \ker \xi_{F_\infty}(T \otimes \rho) \to X_{\infty,*}(F_\infty, T \otimes \rho) \to a^1_\Lambda(\dot{X}_{\infty,*}(F_\infty, \check{T} \otimes \rho^{-1}))$$
$$\to \mathrm{coker}\, \xi_{F_\infty}(T \otimes \rho)$$

*où*

$$\xi_{F_\infty}(T \otimes \rho) : \mathcal{U}^{G_{F_\infty}}_\rho \to \prod_{v \in S^{nd}(F_\infty)} \mathfrak{Z}_\rho(F_{\infty,v}) .$$

*Si de plus $\rho^{-1}$ est admissible pour $V$ et $F_\infty$, on a la suite exacte*

$$0 \to \ker \xi_{F_\infty}(T \otimes \rho) \to X_{\infty,*}(F_\infty, T \otimes \rho) \to a^1_\Lambda(\dot{X}_{\infty,*}(F_\infty, \check{T} \otimes \rho^{-1}))$$
$$\to \mathrm{coker}\, \xi_{F_\infty}(T \otimes \rho) \to \widehat{\check{\mathcal{U}}^{G_{F_\infty}}_\rho} .$$

*En particulier, $X_{\infty,*}(F_\infty, T \otimes \rho)$ est lui aussi de $\Lambda$-torsion.*

*Démonstration.* Par passage à la limite projective des isomorphismes

$$C_{F_n}(T \otimes \rho) : \widehat{H^1_*(F_n, \check{\mathcal{U}}_\rho)} \to H^1_*(F_n, \mathcal{U}_\rho) ,$$

on obtient un homomorphisme de $\Lambda$-modules

$$X_*(F_\infty, T \otimes \rho) \to \varprojlim_n H^1_*(F_n, \mathcal{U}_\rho) \to a^1_\Lambda(\dot{X}_*(F_\infty, \check{T} \otimes \rho^{-1})) .$$

La première flèche est bijective. Les noyau et conoyau de la seconde sont décrits en 2.6.1 ainsi que dans la remarque qui le suit. $\square$

**3.2.2. Corollaire.** *Soit $F_\infty/F$ une $\mathbb{Z}_p$-extension. On suppose vérifiée $(\mathrm{Tors}(F_\infty, V))$. Les applications $C_{F_\infty}(T \otimes \rho)$, pour un caractère continu admissible $\rho$ de $\mathrm{Gal}(F_\infty/F)$ à valeurs dans $\mathbb{Z}_p^*$, induisent par twist par $\rho^{-1}$ un quasi-isomorphisme indépendant de $\rho$*

$$X_{\infty,f}(F_\infty, T) \xrightarrow{\sim} a^1_\Lambda(\dot{X}_{\infty,f}(F_\infty, \check{T}))$$

*et on a plus précisément la suite exacte*

(3.2.1) $\quad 0 \to \ker \xi_{F_\infty}(T) \to X_{\infty,f}(F_\infty, T) \to a^1_\Lambda(\dot{X}_{\infty,f}(F_\infty, \check{T})) \to \mathrm{coker}\, \xi_{F_\infty}(T) \to \widehat{\check{\mathcal{U}}^{G_{F_\infty}}}$

*où*

$$\xi_{F_\infty}(T) : \mathcal{U}^{G_{F_\infty}} \to \prod_{v \in S^{nd}(F_\infty)} \mathfrak{Z}_\rho(F_{\infty,v})$$

*En particulier, $X_{\infty,f}(F_\infty, T)$ est lui aussi de $\Lambda$-torsion.*

*Démonstration.* Il existe sous ces hypothèses un caractère $\rho$ de $\mathrm{Gal}(F_\infty/F)$ admissible pour $T$ tel que $\rho^{-1}$ soit admissible pour $\check{T}$. L'indépendance par rapport à $\rho$ se déduit de 3.1.1 et du fait que le calcul de l'adjoint d'un module $M$ peut se faire en utilisant uniquement les quotient $\widehat{M^{\Gamma_n}}/p^n\widehat{M^{\Gamma_n}}$. □

3.2.3. **Corollaire** (Greenberg). *Le plus grand sous-$\Lambda$-module fini de $X_{\infty,f}(F_\infty,T)$ est égal au noyau de $\xi_{F_\infty}(T)$. Si $\ker \xi_{F_\infty}(T)$ est nul (par exemple si $\mathcal{U}^{G_F}$ est nul, ce qui est équivalent à la nullité de $\mathcal{U}^{G_{F_\infty}}$), $X_{\infty,f}(F_\infty,T)$ n'a pas de sous-modules finis non nuls et est de dimension projective inférieure ou égale à 1.*

Par exemple, s'il existe une place $v \nmid p$ de $S$ telle que $V^{G_{F_\infty,v}} = 0$, $\xi_{F_\infty}(T)$ est injective et $X_{\infty,f}(F_\infty,T)$ est de dimension projective inférieure ou égale à 1. Si $V$ est la représentation $p$-adique associée à une courbe elliptique, cela est le cas s'il existe une place $v \nmid p$ où $E$ a mauvaise réduction additive. On retrouve le résultat démontré par Greenberg ([12, Proposition 4.15])

Il est commode de travailler avec l'application $\Lambda$-sesquilinéaire qui se déduit de $\mathcal{C}_{F_\infty}(T \otimes \rho)$ :

$$\mathrm{Cassels}_{F_\infty}(T \otimes \rho) : X_{\infty,*}(F_\infty, T \otimes \rho) \times X_{\infty,*}(F_\infty, \check{T} \otimes \rho^{-1})) \to \Lambda$$

ou

$$\mathrm{Cassels}_{F_\infty}(T) : X_{\infty,f}(F_\infty, T) \times X_{\infty,f}(F_\infty, \check{T})) \to \Lambda \ .$$

On a donc pour l'une ou l'autre

$$\mathrm{Cassels}(\lambda x, y) = \mathrm{Cassels}(x, \dot{\lambda} y) = \lambda \,\mathrm{Cassels}(x, y) \ .$$

3.2.4. **Proposition.** *On a*

$$\mathrm{Cassels}_{F_\infty}(T)(x,y) = \mathrm{Cassels}_{F_\infty}(\check{T})(y,x)\dot{} \ .$$

*Démonstration.* Il suffit de démontrer l'égalité

$$\mathrm{Cassels}_{F_\infty}(T \otimes \rho)(x,y) = \mathrm{Cassels}_{F_\infty}(\check{T} \otimes \rho^{-1})(y,x)\dot{}$$

pour $\rho$ un caractère de $\mathrm{Gal}(F_\infty/F)$ admissible. Posons $M = X_*(F_\infty, T \otimes \rho)$ et $M' = \dot{X}_*(F_\infty, \check{T} \otimes \rho^{-1})$. L'application $C_{F_\infty}(T \otimes \rho)$ est définie par passage à la limite des $\Lambda$-homomorphismes

$$\widehat{H^1_*(F_n, \check{\mathcal{U}}_\rho)} \to H^1_*(F_n, \mathcal{U}_\rho)$$

et on a le diagramme commutatif

$$\begin{array}{ccc}
\widehat{H^1_*(F_n,\check{\mathcal{U}}_\rho)} & \stackrel{C_{F_n}(T\otimes\rho)}{\to} & H^1_*(F_n,\mathcal{U}_\rho) \\
\uparrow & & \downarrow \\
M_{\Gamma_n} & \to & \widehat{\dot{M}'}^{\Gamma_n} \\
\| & & \uparrow \\
M_{\Gamma_n} & \to & a^1_\Lambda(M')_{\Gamma_n}
\end{array}$$

En prenant le dual de Pontryagin de ce diagramme, on obtient le diagramme commutatif

$$\begin{array}{ccc}
\widehat{H^1_*(F_n,\check{\mathcal{U}}_\rho)} & \stackrel{\dot{C}_{F_n}(\check{T}\otimes\rho^{-1})}{\to} & H^1_*(F_n,\mathcal{U}_\rho) \\
\uparrow & & \downarrow \\
M'_{\Gamma_n} & \to & \widehat{\dot{M}}^{\Gamma_n} \\
\downarrow & & \| \\
\widehat{a^1_\Lambda(\dot{M}')_{\Gamma_n}} & \to & \widehat{M_{\Gamma_n}}
\end{array}$$

On passe ensuite à la limite projective. Les applications $M'_{\Gamma_n} \to a_\Lambda(\widehat{\dot{M}'})_{\Gamma_n}$ induisent alors l'homomorphisme naturel $M' \to a_\Lambda(a_\Lambda(M'))$ et les $C_{F_n}(\check{T} \otimes \rho^{-1})$ induisent l'application $C_{F_\infty}(\check{T} \otimes \rho^{-1})$. □

3.3. **Dualité : cas d'une $\mathbb{Z}_p^2$-extension.** Soit $F_\infty/F$ une $\mathbb{Z}_p^2$-extension. On suppose $(\mathrm{Hyp}(F_\infty, V))$. Soit $L_\infty$ une sous-$\mathbb{Z}_p$-extension de $F_\infty/F$. Disons que $L_\infty$ est admissible si $X_{\infty,f}(L_{\infty,n}, \check{T})$ est un $\Lambda_{L_\infty,n}$-module de torsion pour tout entier $n$. Une telle $\mathbb{Z}_p$-extension existe. En effet, cela revient à montrer que si $F$ est un élément de $\mathbb{Z}_p[[T_1, T_2]]$ (en l'occurrence la série caractéristique de $X_{\infty,f}(F_\infty \check{T})$), il existe un entier $b$ tel que $F(\zeta_n(1+T_2)^b - 1, T_2) \neq 0$ pour tout entier $n$. Dans le cas contraire, $F(T_1, T_2)$ serait divisible par $1 + T_1 - \zeta_{n(b)}(1+T_2)^b$ pour tout entier $b$ avec $n(b)$ entier dépendant de $b$. Ces éléments étant premiers entre eux, cela impliquerait que $F(T_1, T_2)$ est identiquement nul.

On peut alors appliquer le corollaire 3.2.2 à $L_{\infty,n}/L'_n$ : il existe une famille de $\Lambda_{L_\infty,n}$-homomorphismes

$$\mathcal{C}_{L_{\infty,n}}(T) : X_{\infty,f}(L_{\infty,n}, T) \to a^1_{\Lambda_{L_\infty,n}}(\dot{X}_{\infty,f}(L_{\infty,n}, \check{T}))$$

puis définir par passage à la limite projective et en composant avec $r_\infty$ (2.6.4) un homomorphisme de $\Lambda$-modules

$$X_{\infty,f}(F_\infty, T) \to a^1_\Lambda(\dot{X}_{\infty,f}(F_\infty, \check{T}))$$

3.3.1. **Proposition.** *On suppose vérifiés* $(\mathrm{Hyp}(F_\infty, T))$ *et* $(\mathrm{Tors}(F_\infty, V))$. *Le $\Lambda$-homomorphisme $\mathcal{C}_{F_\infty}(T)$ se trouve dans une suite quasi-exacte*

$$0 \to X_{\infty,f}(F_\infty, T) \to a^1_\Lambda(\dot{X}_{\infty,f}(F_\infty, \check{T})) \to \prod_{w \in S^{nd}(F_\infty)} \mathfrak{Z}(F_{\infty,w}) \to 0$$

*Le noyau de $X_{\infty,f}(F_\infty, T) \to a^1_\Lambda(\dot{X}_{\infty,f}(F_\infty, \check{T}))$ est contenu dans $\widehat{\mathcal{U}^{G_{F_\infty}}}$.*

*Démonstration.* Considérons les suites exactes (3.2.1) relatives à la $\mathbb{Z}_p$-extension $L_{\infty,n}/L'_n$ et passons à la limite projective. Les $\ker \xi_{L_{\infty,n}}(T)$ sont finis et d'ordre borné et leur limite projective est finie. La limite projective des $X_{\infty,f}(L_{\infty,n}, \check{T})$ est $X_{\infty,f}(F_\infty, \check{T})$. La limite projective des $a^1_{\Lambda_{L_\infty,n}}(\dot{X}_{\infty,f}(L_{\infty,n}, T))$ est étudiée dans la proposition 2.6.4 (on peut utiliser la remarque qui suit car le plus grand module fini de $X_{\infty,f}(L_{\infty,n}, \check{T})$ est $\ker \xi_{L_\infty,n}$ qui est d'ordre borné par rapport à $n$) : on a la suite quasi-exacte

$$0 \to \varprojlim_n a^1_{\Lambda_{L_\infty,n}}(X_{\infty,*}(L_{\infty,n}, \check{T})) \to a^1_\Lambda(X_{\infty,f}(F_\infty, \check{T}))$$

$$\to \prod_{v \in S^{nd}(F_\infty/L_\infty)} \mathfrak{Z}(F_{\infty,v}) \to 0$$

où $S^{nd}(F_\infty/L_\infty)$ désigne les places de $F_\infty$ non totalement décomposées dans $F_\infty/L_\infty$. Comme $\mathcal{U}^{G_{F_\infty}}$ est supposé fini, la limite projective $W$ du quatrième terme de la suite exacte est quasi-isomorphe à la limite projective des $\prod_{v \in S^{nd}(L_\infty,n)} \mathfrak{Z}(L_{\infty,n,v})$. Soit $v \in S^{nd}(L_{\infty,n})$ ne divisant pas $p$. Elle n'est pas totalement décomposée dans $L_{\infty,n}$, elle est donc totalement décomposée dans $F_\infty/L_\infty$. Comme d'autre part $\mathfrak{Z}(F_{\infty,w})$ est fini, les groupes $\mathfrak{Z}(L_{\infty,n,w})$ sont stationnaires pour $n \gg 0$ et l'application de corestriction de $\prod_{w \in S(L_\infty,n), w|v} \mathfrak{Z}(L_{\infty,n,w})$ est surjective. La limite projective est $\prod_{w \in S(F_\infty), w|v} \mathfrak{Z}(F_{\infty,w})$. On en déduit que $W$ est quasi-isomorphe à $\prod_{w \in S(F_\infty), w|v \in S^{nd}(L_\infty)} \mathfrak{Z}(L_{\infty,n,w})$. Enfin, la limite projective des $\widehat{\mathcal{U}^{G_{L_\infty,n}}}$ est finie. $\square$

3.3.2. **Corollaire.** *On suppose* $(\mathrm{Hyp}(F_\infty, V))$ *et* $(\mathrm{Tors}(F_\infty, V))$. *Le $\Lambda$-module $X_{\infty,f}(F_\infty, T)$ vérifie la propriété (A) : il n'a pas de sous-modules pseudo-nuls non finis et les $a^i_\Lambda(X_{\infty,f}(F_\infty, T))$ sont finis pour $i \geq 2$.*

## 4. Conséquences

Nous pouvons maintenant appliquer les résultats de §1.

4.1. **Descente : $\mathbb{Z}_p$-extension.** Soit $F_\infty/F$ une $\mathbb{Z}_p$-extension. Réécrivons le diagramme (1.2.1) pour $M = X_{\infty,f}(F_\infty, T)$, $M' = \dot{X}_{\infty,f}(F_\infty, \check{T})$ et pour le $\Lambda$-homomorphisme

$$X_{\infty,f}(F_\infty, T) \to a^1_\Lambda(\dot{X}_{\infty,f}(F_\infty, \check{T}))$$

que l'on a construit dans les §3.2 et 3.3. Posons $S_p(T) = H^1_f(F, \mathcal{U})$, $S_p(\check{T}) = H^1_f(F, \check{\mathcal{U}})$, $\check{S}_p(T) = H^1_f(F, T)$ et $\check{S}_p(\check{T}) = H^1_f(F, \check{T})$.

$$\begin{array}{ccccccc}
& & & & \check{S}_p(\check{T}) & & \\
& & & & \downarrow & & \\
0 & \to & S_p(\check{T})/\operatorname{div} & & \operatorname{Hom}_{\mathbb{Z}_p}(M'_\Gamma, \mathbb{Z}_p) & \to & 0 \\
& & \downarrow & & \downarrow & & \\
0 & \to & \widehat{t_{\mathbb{Z}_p}(M'_\Gamma)} & \to & a_\Lambda^1(M')_\Gamma & \to & \operatorname{Hom}_{\mathbb{Z}_p}(M'^\Gamma, \mathbb{Z}_p) & \to & 0 \\
& & \uparrow & & \uparrow & & \uparrow & & \\
0 & \to & t_{\mathbb{Z}_p}(M_\Gamma) & \to & M_\Gamma & \to & L_{\mathbb{Z}_p}(M_\Gamma) & \to & 0 \\
& & \downarrow & & \downarrow & & \downarrow & & \\
0 & \to & \widehat{S_p(T)/\operatorname{div}} & \to & \widehat{S_p(T)} & \to & \operatorname{Hom}_{\mathbb{Z}_p}(\check{S}(T), \mathbb{Z}_p) & \to & 0
\end{array}$$

On peut démontrer en suivant les flèches que l'on retrouve l'application de Cassels, autrement dit que l'on a le diagramme commutatif

$$C_F(T) \quad \begin{array}{ccc} S_p(\check{T})/\operatorname{div} & \to & \widehat{t_{\mathbb{Z}_p}(M'_\Gamma)} \\ \uparrow & & \uparrow \\ \widehat{S_p(E)/\operatorname{div}} & \leftarrow & t_{\mathbb{Z}_p}(M_\Gamma) \end{array}$$

D'autre part, lorsque $\check{S}(T)$ n'est pas fini, on obtient une forme bilinéaire $\langle \cdot, \cdot \rangle_\gamma$

$$\check{S}(T) \times \check{S}(\check{T}) \to \mathbb{Q}_p \ .$$

Elle dépend de $\gamma$ et bien sûr de la $\mathbb{Z}_p$-extension $F_\infty$. Pour $\rho$ caractère non trivial de $\Gamma$ à valeurs dans $\mathbb{Z}_p^*$, on note $\langle \cdot, \cdot \rangle_\rho = (\log_p \rho(\gamma))^{-1} \langle \cdot, \cdot \rangle_\gamma$. On peut démontrer que l'on retrouve la hauteur $p$-adique ordinaire associée à $\rho$ (cf. [23] dans un cadre un peu différent). Nous n'en aurons pas besoin.

4.2. **Descente : $\mathbb{Z}_p^2$-extension.** Maintenant que nous avons construit $\mathcal{C}_{F_\infty}(T)$ en nous aidant de coinvariants convenables (c'est-à-dire de torsion pour la $\mathbb{Z}_p$-extension corespondante), il est possible de redescendre en utilisant les homomorphismes fonctoriels construits dans le §1 et plus particulièrement §1.2. Ce qui permet d'obtenir des informations pour les $\mathbb{Z}_p$-extensions telles que $X_{\infty,f}(L_\infty, T)$ n'est pas de torsion.

On fait les hypothèses (Hyp($F_\infty, V$)), (Hyp($F_\infty, \check{V}$)) et (Tors($F_\infty, V$)). Soit $\mathfrak{z}(T) = \mathfrak{z}(F_\infty, T) = \prod_{w \in S^{nd}(F_\infty)} \mathfrak{z}(F_{\infty,w}, T)$ (pour le caractère $\rho$ trivial). On a contruit dans le §3.3 la suite quasi-exacte de $\Lambda$-modules suivante

$$0 \to X_{\infty,f}(F_\infty, T) \to a_\Lambda^1(\dot{X}_{\infty,f}(F_\infty, \check{T})) \to \mathfrak{z}(F_\infty, T) \to 0$$

Soit $L_\infty$ une sous-$\mathbb{Z}_p$-extension de $F_\infty/F$. Posons $\Gamma' = \operatorname{Gal}(F_\infty/L_\infty)$, $\Lambda_{L_\infty} = \mathbb{Z}_p[[\operatorname{Gal}(L_\infty/F)]]$. Dans le cas où $X_{\infty,f}(L_\infty, T)$ est de torsion, on a alors le diagramme commutatif suivant dont les lignes et les colonnes sont quasi-exactes

$$\begin{array}{ccccccccc}
& & 0 & & 0 & & & & \\
& & \downarrow & & \uparrow & & & & \\
& & \widehat{\mathfrak{z}(\check{T})_{\Gamma'}} & & a_{\Lambda_{L_\infty}}^1(\widehat{\mathfrak{z}(T)_{\Gamma'}}) & & & & \\
& & \downarrow & & \uparrow & & & & \\
0 & \to & \mathfrak{z}(T)^{\Gamma'} & \to & X_{\infty,f}(F_\infty, T)_{\Gamma'} & \to & a_\Lambda^1(\dot{X}_{\infty,f}(F_\infty, \check{T})_{\Gamma'}) & \to & \mathfrak{z}(T)_{\Gamma'} & \to & 0 \\
& & & & \downarrow & & \uparrow & & & & \\
& & & & X_{\infty,f}(L_\infty, T) & \to & a_\Lambda^1(\dot{X}_{\infty,f}(L_\infty, \check{T})) & & & & \\
& & & & \downarrow & & \uparrow & & & & \\
& & & & 0 & & 0 & & & & 
\end{array}$$

et le quasi-isomorphisme du §3.2

$$X_{\infty,f}(L_\infty, T) \to a_{\Lambda_{L_\infty}}^1(\dot{X}_{\infty,f}(L_\infty, T))$$

Ne supposons plus $X_{\infty,f}(L_\infty,\check{T})$ de $\Lambda_{L_\infty}$-torsion. On a alors le diagramme commutatif suivant dont les lignes sont quasi-exactes :

$$\begin{array}{ccccc}
a^1_{L_\infty}(t_{\Lambda_{L_\infty}}(\dot{X}_{\infty,f}(F_\infty,\check{T})_{\Gamma'})) & & & & \mathrm{Hom}_{\Lambda_{L_\infty}}(\dot{X}_{\infty,f}(F_\infty,\check{T})_{\Gamma'},\Lambda_{L_\infty}) \\
\uparrow\sim & & & & \downarrow \\
0 \to a^1_{L_\infty}(\dot{X}_{\infty,f}(F_\infty,\check{T})) & \to & a^1_{L_\infty}(\dot{X}_{\infty,f}(F_\infty,\check{T})_{\Gamma'}) & \to & \mathrm{Hom}_{\Lambda_{L_\infty}}(\dot{X}_{\infty,f}(F_\infty,\check{T})^{\Gamma'},\Lambda_{L_\infty}) \to 0 \\
\uparrow & & \uparrow & & \uparrow \\
0 \to t_{\Lambda_{L_\infty}}(X_{\infty,f}(F_\infty,T)_{\Gamma'}) & \to & X_{\infty,f}(F_\infty,T)_{\Gamma'} & \to & L_{\Lambda_{L_\infty}}(X_{\infty,f}(F_\infty,T)_{\Gamma'}) \to 0
\end{array}$$

On en déduit alors comme du §1.2 des homomorphismes

$$t_{\Lambda_{L_\infty}}(X_{\infty,f}(F_\infty,T)_{\Gamma'}) \to \dot{a}_{L_\infty}(X_{\infty,f}(F_\infty,\check{T})_{\Gamma'})$$

et une forme sesqui-linéaire

$$\mathrm{Cassels}_{L_\infty}(T) : t_{\Lambda_{L_\infty}}(X_{\infty,f}(F_\infty,T)_{\Gamma'}) \times t_{\Lambda_{L_\infty}}(X_{\infty,f}(F_\infty,\check{T})_{\Gamma'}) \to \mathrm{Frac}(\Lambda_{L_\infty})/\Lambda_{L_\infty}$$

quasi non dégénérée vérifiant

$$\mathrm{Cassels}_{L_\infty}(T)(\gamma x, y) = \gamma\,\mathrm{Cassels}_{L_\infty}(T)(x,y) = \mathrm{Cassels}_{L_\infty}(T)(x,\gamma^{-1}y)$$
$$\mathrm{Cassels}_{L_\infty}(T)(x,y) = \mathrm{Cassels}_{L_\infty}(\check{T})(y,x)\ .$$

On obtient aussi une hauteur $p$-adique qui est un accouplement sur les quotients sans $\Lambda_{L_\infty}$-torsion de $X_{\infty,f}(F_\infty,T)_{\Gamma'}$ et de $X_{\infty,f}(F_\infty,\check{T})_{\Gamma'}$ à valeurs dans $\Lambda_{L_\infty}$.

Le noyau (resp. conoyau) de $X_{\infty,f}(F_\infty,T)_{\Gamma'} \to a^1_{L_\infty}(\dot{X}_{\infty,f}(F_\infty,\check{T})_{\Gamma'})$ est de torsion et quasi-isomorphe à $\mathfrak{Z}(T)^{\Gamma'}$(resp. $\mathfrak{Z}(T)_{\Gamma'}$) qui est d'ailleurs annulé par une puissance de $p$. On en déduit une suite exacte

$$0 \to \mathfrak{Z}(T)^{\Gamma'} \to t_{\Lambda_{L_\infty}}(X_{\infty,f}(F_\infty,T)_{\Gamma'}) \to a^1_{L_\infty}(\dot{X}_{\infty,f}(F_\infty,\check{T})_{\Gamma'}) \to Z \to 0$$

avec $Z$ annulé par une puissance de $p$ et de $\mu$-invariant inférieur à celui de $\mathfrak{Z}(T)^{\Gamma'}$. La série caractéristique de $t_{\Lambda_{L_\infty}}(X_{\infty,f}(F_\infty,\check{T})_{\Gamma'})$ divise donc celle de $t_{\Lambda_{L_\infty}}(X_{\infty,f}(F_\infty,T)_{\Gamma'})$. Par symétrie, on en déduit qu'elles sont égales et que l'on a la suite quasi-exacte

$$0 \to \mathfrak{Z}(T)^{\Gamma'} \to t_{\Lambda_{L_\infty}}(X_{\infty,f}(F_\infty,T)_{\Gamma'}) \to a^1_{L_\infty}(\dot{X}_{\infty,f}(F_\infty,\check{T})_{\Gamma'}) \to \mathfrak{Z}(T)_{\Gamma'} \to 0$$

Comme $\widehat{\mathfrak{Z}(\check{T})_{\Gamma'}}$ et $\widehat{\mathfrak{Z}(T)_{\Gamma'}}$ sont de torsion, les suites suivantes sont quasi-exactes

$$0 \to \widehat{\mathfrak{Z}(\check{T})_{\Gamma'}} \to t_{\Lambda_{L_\infty}}(X_{\infty,f}(F_\infty,T)_{\Gamma'}) \to t_{\Lambda_{L_\infty}}(X_{\infty,f}(L_\infty,T)) \to 0$$
$$0 \to a^1_{\Lambda_{L_\infty}}(t_{\Lambda_{L_\infty}}(\dot{X}_{\infty,f}(L_\infty,\check{T}))) \to a^1_{\Lambda_{L_\infty}}(t_{\Lambda_{L_\infty}}(\dot{X}_{\infty,f}(F_\infty,\check{T})_{\Gamma'})) \to a^1_{\Lambda_{L_\infty}}(\widehat{\mathfrak{Z}(T)_{\Gamma'}}) \to 0$$

En remarquant que $\mathfrak{Z}(T)^{\Gamma'}$ et $\widehat{\mathfrak{Z}(\check{T})_{\Gamma'}}$ ont même série caractéristique (cela peut se voir soit sur le diagramme, soit directement : seuls les nombres de Tamagawa aux places de $S$ totalement décomposées dans $L_\infty$ interviennent et on a alors $\mathrm{Tam}_v(T) = \mathrm{Tam}_v(\check{T})$ pour une place ne divisant pas $p$, en fait $\mathfrak{Z}(T)$ et $\mathfrak{Z}(\check{T})$ sont quasi-isomorphes), on en déduit le théorème :

**4.2.1. Théorème.** *Supposons* $(\mathrm{Hyp}(F_\infty,V))$, $(\mathrm{Hyp}(F_\infty,\check{V}))$ *et* $(\mathrm{Tors}(F_\infty,V))$. *Soit* $f(F_\infty,T)$ *la série caractéristique de* $X_{\infty,f}(F_\infty,T)$ *et* $f(F_\infty,\check{T})$ *la série caractéristique de* $X_{\infty,f}(F_\infty,\check{T})$. *Alors*

$$\dot{f}(F_\infty,\check{T})\Lambda = f(F_\infty,T)\Lambda$$

*Pour toute sous-$\mathbb{Z}_p$-extension de $F_\infty/F$, soit $f^*(L_\infty,\check{T})$ (resp. $f^*(L_\infty,T)$) la série caractéristique du sous-module de torsion de $X_{\infty,f}(L_\infty,\check{T})$ (resp. $X_{\infty,f}(L_\infty,T)$). Alors*

$$\dot{f}^*(L_\infty,\check{T})\Lambda_{L_\infty} = f^*(L_\infty,T)\Lambda_{L_\infty}$$

*Autrement dit, pour tout caractère $\rho$ de $\mathrm{Gal}(L_\infty/F)$ à valeurs dans $\mathbb{C}_p^*$, on a*

$$\rho(f^*(L_\infty,T)) = \rho^{-1}(f^*(L_\infty,\check{T}))$$

## 5. La situation diédrale

**5.1. Préliminaires.** Soit $K$ un corps quadratique imaginaire de discriminant $D$ et $p$ un nombre premier impair et premier à $D$. Il existe une unique extension $K_\infty$ de $K$ dont le groupe de Galois est topologiquement isomorphe à $\mathbb{Z}_p^2$. Elle contient deux $\mathbb{Z}_p$-extensions qui sont la sous-$\mathbb{Z}_p$-extension cyclotomique $K\mathbb{Q}_\infty$ de $K$ et la sous-extension anticyclotomique (diédrale sur $\mathbb{Q}$) $H_\infty = D_\infty$ de $K[p^\infty] = \cup K[p^n]$, le Ringklasskörper de K de rayon une puissance de $p$. Soit $G_\infty = \mathrm{Gal}(K_\infty/K)$.

L'algèbre d'Iwasawa associée est $\Lambda = \Lambda_{K_\infty} = \mathbb{Z}_p[[G_\infty]]$. On a une dualité

$$\Lambda_{K_\infty} \times \mathrm{Hom}(G_\infty, \mathbb{C}_p^\times) \to \mathbb{C}_p \ .$$

Ce qui permet de voir les éléments de $\Lambda_{K_\infty}$ comme des fonctions sur $\mathrm{Hom}(G_\infty, \mathbb{C}_p^\times)$ à valeurs dans $\mathbb{C}_p$. Tout élément de $\mathrm{Hom}(G_\infty, \mathbb{Z}_p^\times)$ est de la forme $\nu^a \chi^b$ avec $\chi$ la $p$-partie du caractère cyclotomique et $\nu$ un caractère diédral. On note $\chi_{cycl}$ le caractère cyclotomique et $\nu_{died}$ un caractère diédral.

La théorie d'Iwasawa d'une courbe elliptique sur un corps quadratique imaginaire et des famille des points de Heegner tire ses origines de l'article de Mazur ([20], voir aussi Kurcanov, [**?**]). Grâce aux résultats récents de Vatsal et Cornut, un regain d'intérêt s'est manifesté. Mais il y a bien d'autres résultats montrés ou en voie de l'être et je voudrais les placer ici un peu plus dans le contexte de la théorie d'Iwasawa.

**5.2. Théorie arithmétique.** Soit $E$ une courbe elliptique définie sur $\mathbb{Q}$ ou $K$ de conducteur $N$ ayant bonne réduction ordinaire en $p$. Soit $T = T_p(E)$ son module de Tate, c'est-à-dire la limite projective des points de $p^n$-torsion pour $n$ entier et $V_p(E) = \mathbb{Q}_p \otimes T_p(E)$. Soit $L$ une extension finie de $K$. Le groupe de Selmer $S_p(E/L)$ de $E/L$ vérifie la suite exacte

$$0 \to \mathbb{Q}_p/\mathbb{Z}_p \otimes E(L) \to S_p(E/L) \to \text{III}(E/L)(p) \to 0$$

Son dual de Pontryagin $\hat{S}_p(E/L)$ est $\mathrm{Hom}_{\mathbb{Z}_p}(S_p(E/L), \mathbb{Q}_p/\mathbb{Z}_p)$. Sa variante compacte (voir [2], [10]) $\check{S}_p(E/L)$ est la limite projective des groupes de Selmer relatif à la multiplication par $p^n$ :

$$0 \to \mathbb{Z}_p \otimes_{\mathbb{Z}} E(L) \to \check{S}_p(E/L) \to T_p(\text{III}(E/L)) \to 0$$

Lorsque $\text{III}(E/L)(p)$ est fini, le dernier terme est nul. Avec les notations du §2, on a $\check{S}_p(E/L) = H^1_f(L, T_p(E))$ et $S_p(E/L) = H^1_f(L, V_p(E)/T_p(E))$. Le quotient de $\text{III}(E/L)(p)$ par sa partie divisible est le groupe de Tate-Shafarevich $\text{III}(T_p(E)/L)$ associé à la représentation $p$-adique $T_p(E)$ et vaut $S_p(E/L)/\mathbb{Q}_p/\mathbb{Z}_p \otimes \check{S}_p(E/L)$. Ainsi,

$$S_p(E/K_\infty) = \varinjlim_{L} S_p(E/L) = H^1_f(K_\infty, \mathcal{U})$$
$$\widehat{S_p(E/K_\infty)} = X_{\infty,f}(K_\infty, T_p(E))$$
$$\check{S}_p(E/K_\infty) = \varprojlim_{L} \check{S}_p(E/L)$$

où la limite projective est prise relativement aux applications de norme (corestriction). Enfin, il n'est pas difficile de montrer que pour une $\mathbb{Z}_p$-extension, par exemple $D_\infty$, on a

$$\check{S}_p(E/D_\infty) = \mathrm{Hom}_{\Lambda_{D_\infty}}(\widehat{S_p(E/D_\infty)}, \Lambda_{D_\infty}) \cong \mathrm{Hom}_{\Lambda_{D_\infty}}(\widehat{S_p(E/D_\infty)}/t_{\Lambda_{D_\infty}}(\widehat{S_p(E/D_\infty)}), \Lambda_{D_\infty})$$

et que

$$\widehat{S_p(E/D_\infty)}/t_{\Lambda_{D_\infty}}(\check{S}_p(E/D_\infty)) \to \mathrm{Hom}(\widehat{S_p(E/D_\infty)}, \Lambda_{D_\infty})$$

est injectif avec conoyau fini. En particulier, $\check{S}_p(E/D_\infty)$ est le $\Lambda_{D_\infty}$-dual de $X_{\infty,f}(D_\infty, T_p(E)) = \widehat{S_p(E/D_\infty)}$ et est libre.

**5.2.1. Théorème** (Kato). *Si $E$ est définie sur $\mathbb{Q}$, $\check{S}_p(E/K\mathbb{Q}_\infty)$ est de torsion sur $\Lambda_{\mathbb{Q}_\infty K}$ et $\check{S}_p(E/K_\infty)$ est de torsion sur $\Lambda_{K_\infty}$. Il en est de même de $\widehat{S_p(E/K_\infty)}$ .*

Il suffit d'appliquer le théorème démontré par Kato dans [18] à $E$ et à sa tordue par le caractère quadratique de $K/\mathbb{Q}$.

Soit $\mathcal{L}_p(E/K_\infty)$ une série caractéristique du module de torsion $\widehat{S_p(K_\infty/K)}$.

5.3. **Théorie analytique.** Soit $L_p(E/K_\infty)$ la fonction $L$ $p$-adique interpolant les valeurs $L(E,\rho,1)$ pour $\rho$ caractère d'ordre fini de $\mathrm{Gal}(K_\infty/K)$. On peut trouver sa définition dans [24] qui suit de très près une construction antérieure de Hida. D'autres constructions ont été faites par Bertolini et Darmon [5]. Par un théorème de Rohrlich [27], $L_p(E/K_\infty)$ est non nulle.

**Conjecture** (Conjecture principale, [25]). *Les idéaux de $\Lambda_{K_\infty}$ engendrés par $L_p(E/K_\infty)$ et par $\mathcal{L}_p(E/K_\infty))$ sont égaux.*

**Remarque.** On peut utiliser le théorème de Kato pour obtenir une divisibilité lorsqu'on se restreint à $\mathrm{Gal}(K\mathbb{Q}_\infty/K)$.

5.4. **Equation fonctionnelle.** Soit $c$ une conjugaison complexe induisant l'automorphisme non trivial de $\mathrm{Gal}(K/\mathbb{Q})$. Elle agit sur le groupe des caractères $\hat{G}_\infty$ de $G_\infty$ : $\rho^c(\tau) = \rho(c\tau c^{-1})$. Ainsi, si $\chi$ se factorise par $\mathrm{Gal}(K\mathbb{Q}_\infty/K)$, on a $\chi^c = \chi$. Si $\nu$ est diédral, $\nu^c = \nu^{-1}$. On considère l'involution suivante sur $\hat{G}_\infty$ : $\rho^\iota = \rho^{-c}$. Ainsi, $\rho^\iota(\tau) = \rho(c^{-1}\tau c)^{-1}$, $\chi^\iota_{cycl} = \chi^{-1}_{cycl}$, $\nu^\iota_{died} = \nu_{died}$.

Les deux fonctions vérifient une équation fonctionnelle
$$(L(\rho)) = (L(\rho^\iota))$$
Pour la première, cela se déduit de l'équation fonctionnelle complexe et on a en fait
$$L_p(E/K_\infty)(\rho^\iota) = \epsilon_D(-N) L_p(E/K_\infty)(\rho) \ .$$
En appliquant l'automorphisme non trivial $c$ de $K/\mathbb{Q}$ qui laisse stable $E$ ainsi que tous les modules définis et la proposition 4.2.1 on obtient la seconde équation fonctionnelle. On en déduit que l'on peut définir le signe de l'équation fonctionnelle de $\mathcal{L}_p(E/K_\infty)$ : le groupe de cohomologie $H^1(\{1,\iota\}, \Lambda^\times_{K_\infty})$ est d'ordre 2 et admet $-1$ comme élément non trivial. Ainsi, on peut choisir $\mathcal{L}_p(E/K_\infty)$ (qui est alors défini à une unité près de $\mathbb{Z}_p^*$) de manière à ce que
$$\mathcal{L}_p(E/K_\infty)(\rho^\iota) = \epsilon_p \mathcal{L}_p(E/K_\infty)(\rho)$$
avec $\epsilon_p = \pm 1$.

5.4.1. **Proposition.** *Soit $\lambda_0(D_\infty)$ le rang du $\Lambda_{D_\infty}$-module $S(\widehat{E/D_\infty})$. Alors,*
$$\epsilon_p = (-1)^{\mathrm{rg}_{\mathbb{Z}_p} \check{S}_p(E/K)} = (-1)^{\lambda_0(D_\infty)}$$

*Démonstration.* On utilise le théorème de contrôle 2.5.1 et l'argument de Guo et Greenberg que l'on va rappeler. En se restreignant aux caractères cyclotomiques, l'équation fonctionnelle devient
$$\mathcal{L}_p(E/K_\infty)(\chi^{-s}_{cycl}) = \epsilon_p \mathcal{L}_p(E/K_\infty)(\chi^s_{cycl})$$
On obtient de même (Guo et Greenberg)
$$\epsilon_p = (-1)^{\mathrm{rg}_{\mathbb{Z}_p} \check{S}_p(E/K)}$$
□

Expliquons l'argument de Guo [14] tel qu'il a été repris par Greenberg dans [12]. Il s'appuie sur les faits suivants. Soit $\Lambda = \mathbb{Z}_p[[\Gamma]]$, $\Gamma'$ un sous-$\mathbb{Z}_p$-module de $\Gamma$ isomorphe à $\mathbb{Z}_p$. On note $\Xi_n = \Gamma/\Gamma'_n$, $\Lambda_n = \mathbb{Z}_p[[\Xi_n]]$, $\Xi$ un sous-$\mathbb{Z}_p$-module de $\Gamma$ tel que $\Xi \cap \Gamma' = \{0\}$, $\Lambda_\Xi = \mathbb{Z}_p[[\Xi]]$. Soit $M$ un $\Lambda$-module de torsion et de type fini ; $M_{\Gamma'_n}$ est un $\Lambda_\Xi$-module de type fini. Soit $\lambda_n$ le $\Lambda_\Xi$-rang de $M_{\Gamma'_n}$. Alors $\lambda_n$ est stationnaire et on a $\lambda_n \equiv \lambda_0 \bmod p-1$. En effet, les représentations irréductibles de $\Gamma'/\Gamma'_n$ sont de degré divisible par $(p-1)p^{n-1}$.

Soit $\mathfrak{L}_n$ le quotient de $M_n$ par son $\Lambda_\Xi$-module de torsion $\mathfrak{T}_n$ et $\mathfrak{T}_\infty$ le noyau de $M \to \mathfrak{L}_\infty$ où $\mathfrak{L}_\infty$ est la limite projective des $\mathfrak{L}_n$. Alors $\mathfrak{L}_\infty$ est un $\Lambda_\Xi$-module de torsion, l'application naturelle $\mathfrak{L}_\infty \to \mathfrak{L}_n$ est injective pour $n$ assez grand, les rangs de $\mathfrak{L}_\infty$ et de $\mathfrak{L}_n$ sur $\Lambda_\Xi$ sont égaux pour $n$ assez grand et on a donc
$$\mathrm{rg}_{\Lambda_\Xi} \mathfrak{L}_\infty \equiv \lambda_0 \bmod p-1 \ .$$
Supposons qu'il existe une forme bilinéaire alternée sur $\mathfrak{T}_n$ quasi non dégénérée (ou telle que les noyaux soient annulés par une puissance de $p$). Alors le rang de $\mathfrak{T}_\infty$ en tant que $\Lambda_\Xi$-module est pair. Pour le démontrer, on remarque, en utilisant le lemme du serpent et le fait que $\mathfrak{L}_\infty \to \mathfrak{L}_n$ est injective, que l'application $\mathfrak{T}_\infty \to \mathfrak{T}_n$ est surjective. Si $\mathfrak{p}$ est un idéal de $\Lambda$ de hauteur 1 premier à $p$, la forme bilinéaire alternée sur $\mathfrak{T}_n$ induit par localisation une forme alternée non dégénérée sur $\Lambda_{\mathfrak{p}} \otimes \mathfrak{T}_n$. Le $\Lambda_{\mathfrak{p}}$-rang de $\Lambda_{\mathfrak{p}} \otimes \mathfrak{T}_\infty$ est alors pair. En effet, il suffit d'appliquer le lemme suivant :

**Lemme** (Guo). *Soit $A$ un anneau de valuation discrète d'uniformisante $\pi$ et soit $M_n$ un système projectif de $A$-modules de longueur finie tel que $M_\infty = \varprojlim_n M_n$ soit de type fini et que l'application naturelle $M_\infty \to M_n$ soit surjective. Alors, il existe un entier $d$ et des entiers $r_1(n), \cdots, r_d(n)$ tels que $M_n \cong A/\pi^{r_1(n)} \times \cdots \times A/\pi^{r_d(n)}$ avec $r_1(n) \geq \cdots \geq r_d(n)$ et le rang sur $\mathbb{Z}_p$ de $M_\infty$ est égal au nombre d'entiers $j$ tels que la suite $r_j(n)$ soit non bornée.*

Si maintenant $X_n$ est muni d'une forme bilinéaire alternée pour tout entier $n$, on a $r_{2j-1}(n) = r_{2j}(n)$ pour $j \geq 1$ et le rang sur $\mathbb{Z}_p$ de $X_\infty$ est donc pair.

5.5. **La conjecture de Mazur.** On suppose toujours que le discriminant de $K$ est premier au conducteur $N_E$ de $E$. Dans [20], Mazur conjecture que $S_p(\widehat{E/D_\infty})$ est un $\Lambda_{D_\infty}$-module de rang 1 si $\epsilon(-N_E) = -1$ et de rang 0 si $\epsilon(-N_E) = 1$.

Cette conjecture et la proposition 5.4.1 impliquent que les signes des équations fonctionnelles de $\mathcal{L}_p(E/K_\infty)$ et de $L_p(E/K_\infty)$ sont égaux.

Récemment, la situation de cette conjecture a énormément évolué dans le cas de l'hypothèse de Heegner mais aussi dans les autres cas. Nous appellerons hypothèses techniques des hypothèses qui devraient pouvoir être affaiblies ou évoluer rapidement jusqu'à disparaître et que l'on trouvera dans les articles originaux :

5.5.1. **Théorème** (Bertolini-Darmon+Vatsal, [6], [28]). *Supposons que $\epsilon(-N_E) = 1$. Supposons de plus que $\ell^2 \nmid N_E$ si $\ell$ est inerte dans $K$ + des hypothèses techniques, alors $S_p(\widehat{E/D_\infty})$ est de torsion.*

La démonstration est en deux parties :
- démontrer la non nullité de la fonction $L$ $p$-adique $L_p(E/D_\infty)$ (théorème sur les familles de $L(E/K, \eta)$ pour $\eta$ un caractère diédral d'ordre fini et de conducteur une puissance de $p$)
- démontrer que si $L_p(E/D_\infty)$ est non nul, $S_p(\widehat{E/D_\infty})$ est de torsion

5.5.2. **Théorème** (Cornut-Vatsal + Bertolini-Darmon-Nekovář [4]). *Lorsque $\epsilon(-N_E) = -1$ et que tous les nombres premiers divisant $N_E$ sont décomposés dans $K$, $S_p(\widehat{E/D_\infty})$ est de rang 1.*

L'énoncé complet que l'on attend est montré ou sur le point de l'être :

5.5.3. **Théorème.** *Supposons que $\epsilon(-N_E) = -1$ et que $p^2$ ne divise pas $N_E$. Alors, $S_p(\widehat{E/D_\infty})$ est de rang 1.*

La démonstration comporte plusieurs étapes :
- Construire des points $x_n$ de $E(K_n)$ en utilisant une paramétrisation de $E$ par une courbe modulaire ou une courbe de Shimura et les points de Heegner ou les points spéciaux provenant de la théorie de la multiplication complexe $x(p^n)$ (idée de Gross exploitée par Bertolini et Darmon, [5]). En modifiant légèrement ces points, on obtient des points compatibles pour les applications de trace et donc un élément $z_\infty^{spec}$ de $\tilde{S}_p(K_\infty)$ et un sous-module $\mathcal{H}_\infty$ de $\tilde{S}_p(E/D_\infty)$.
- Montrer que $z_\infty^{spec}$ est non nul (conjecture de Mazur), ce qui se ramène facilement à montrer qu'il existe un entier $n$ tel que $x_n$ est non nul. C'est le rôle des théorèmes de Cornut et Vatsal. On pourrait aussi peut-être utiliser les formules démontrées par Zhang ([29]) généralisant les formules de Gross-Zagier et qui sont du type :

$$L'(E/K, \nu, 1) = C <x_n^{(\nu)}, x_n^{(\nu)}>$$

avec $C$ non nul et utiliser un théorème de non-annulation de la famille des $L'(E/K, \nu, 1)$ pour un caractère $\nu$ diédral d'ordre une puissance de $p$. On en déduit alors facilement que $\mathcal{H}_\infty$ est un module libre de rang 1.
- Utiliser les techniques de Kolyvagin [17] pour démontrer que le quotient $\tilde{S}_p(E/D_\infty)/\mathcal{H}_\infty$ est de torsion et donc que $\tilde{S}_p(E/D_\infty)$ et $S_p(\widehat{E/D_\infty})$ sont de $\Lambda$-rang 1. Ici, c'est la notion de système d'Euler qui est fondamentale. La définition précise des points $x_n$ ne sert pas mais le fait qu'il existe des points $x(np^m)$ pour $n$ sans facteurs carrés définis sur le Ringklasskörper de rayon $np^m$ vérifiant des relations convenables.

5.6. **Quelques remarques supplémentaires.** Les techniques de Kolyvagin ont permis à Bertolini [3] de démontrer le théorème suivant :

5.6.1. **Théorème.** *Soit $I(\mathcal{H}_\infty)$ la série caractéristique du $\Lambda_{K_\infty}$-module de torsion $\tilde{S}_p(E/D_\infty)/\mathcal{H}_\infty$. Alors, sous l'hypothèse de Heegner et sous des hypothèses techniques, $(\gamma - 1)I(\mathcal{H}_\infty)$ annule le sous-$\Lambda$-module de torsion $t_{\Lambda_{D_\infty}}(S_p(\widehat{E/D_\infty}))$ de $S_p(\widehat{E/D_\infty})$. De plus si $\gamma - 1$ ne divise pas $I(\mathcal{H}_\infty)$, $I(\mathcal{H}_\infty)$ annule $t_{\Lambda_{D_\infty}}(S_p(\widehat{E/D_\infty}))$.*

Les hypothèses techniques doivent toutes pouvoir être enlevées à condition de rajouter une puissance de $p$. Et on devrait pouvoir remplacer l'hypothèse de Heegner par l'hypothèse que $\epsilon(-N_E) = -1$ en utilisant les points de Heegner-Shimura. Cela faisait suite à la conjecture suivante. Soit $u = (\sharp\mathcal{O}_K^*)/2$ et $c_E$ la constante de Manin correspondant à la paramétrisation de $E$ par $X_0(N_E)$.

5.6.2. **Conjecture** ([25]). *Les deux éléments $c_E^2 u^2 \mathcal{T}(E/D_\infty)$ et $I(\mathcal{H}_\infty)^2$ engendrent le même idéal de $\Lambda_{D_\infty}$.*

Il semble qu'en poussant à fond la méthode de Kolyvagin, on devrait pouvoir montrer la divisibilité de $I(\mathcal{H}_\infty)^2$ par $c_E^2 u^2 \mathcal{T}(E/D_\infty)$. Remarquons qu'on déduit des résultats du §4.2 que $\mathcal{T}(E/D_\infty)$ est un carré, ce qui est compatible avec la conjecture précédente. En effet, en composant avec l'involution $c$ et en identifiant $T_p(E)$ avec $T_p(E)^*(1)$ par l'accouplement alterné de Weil, on obtient une forme bilinéaire alternée

$$t_{\Lambda_{D_\infty}}(X_{\infty,f}(K_\infty, T_p(E))_{\text{Gal}(K_\infty/D_\infty)}) \times t_{\Lambda_{D_\infty}}(X_{\infty,f}(K_\infty, T_p(E))_{\text{Gal}(K_\infty/D_\infty)}) \to \text{Frac}(\Lambda_{D_\infty})/\Lambda_{D_\infty}$$

En tenant compte des noyaux et de la différence entre $X_{\infty,f}(K_\infty, T_p(E))_{\text{Gal}(K_\infty/D_\infty)}$ et $X_{\infty,f}(D_\infty, T_p(E))$ dont la série caractéristique est une puissance de $p$, on en déduit que $\mathcal{T}(E/D_\infty)$ est un carré.

Une conséquence de la démonstration du résultat de Cornut est la suivante :

5.6.3. **Proposition.** *$I(\mathcal{H}_\infty)$ n'est pas divisible par $p$.*

Supposons que $\epsilon(-N_E) = -1$. Nous avons défini précédemment une forme bilinéaire

$$\langle\langle \cdot, \cdot \rangle\rangle_{\chi_{cycl}} : \check{S}_p(E/D_\infty) \times \check{S}_p(E/D_\infty) \to \Lambda_{D_\infty} .$$

Elle peut s'écrire en termes des hauteurs $p$-adiques classiques de la manière suivante

$$\langle\langle x, y \rangle\rangle_{\chi_{cycl}} = \left( \frac{1}{[D_n : K]} \sum_{\sigma,\tau} <\sigma x_n, \tau y_n>_{\chi_n} \sigma\tau^{-1} \right)_n .$$

Ici, $\chi_n$ est un caractère de $\text{Gal}(K_\infty/D_n)$ dont la restriction à $\text{Gal}(K_\infty/D_\infty)$ est $\chi_{cycl}$. On peut reprendre la démonstration de [25] pour démontrer :

5.6.4. **Théorème.** *Soit $\rho$ un caractère diédral de $\text{Gal}(K_\infty/K)$ à valeurs dans $\mathbb{Z}_p^*$.*

*1) Soit $r_{D_\infty}$ le rang de $\check{S}_p(E/D_\infty)$ en tant que $\Lambda_{D_\infty}$-modules. Alors, $\mathcal{L}(E/K_\infty)(\rho\chi_{cycl}^s)$ a un zéro en $\chi_{cycl}$ de multiplicité supérieure ou égale à $r_{D_\infty}$.*

*2) Ce zéro est d'ordre exactement $r_{D_\infty}$ si et seulement si $\langle\langle \cdot, \cdot \rangle\rangle_{\chi_{cycl}}$ est non dégénérée.*

*3) On a dans ce cas*

$$\lim_{s \to 0} \frac{\mathcal{L}(K_\infty/K)(\rho\chi_{cycl}^s)}{s^{r_{D_\infty}}} \sim disc_{\check{S}_p(E/D_\infty)} \langle\langle \cdot, \cdot \rangle\rangle_{\chi_{cycl}}(\rho)\mathcal{T}_p(E/D_\infty)(\rho)$$

*où $\mathcal{T}_p(E/D_\infty)$ est une série caractéristique du $\Lambda_{D_\infty}$-module de torsion de $S_p(\widehat{E/D_\infty}/K)$.*

5.6.1. Les théorèmes ou conjectures précédentes impliquent que $r_{D_\infty}$ est en fait égal à 1 lorsque $\epsilon(-N_E) = -1$. D'autre part, l'ordre du zéro de $\mathcal{L}(E/K_\infty)(\rho\chi_{cycl}^s)$ en $s = 0$ est impair. Il serait intéressant de montrer qu'il existe un point de Heegner $z_n$ dont la hauteur $p$-adique $<z_n, z_n>_{\chi_n}$ est non nulle. Cela n'est connu que si $E$ est à multiplication complexe ([7]).

Revenons sur le module des points de Heegner. Soit $\mathcal{H}_n$ le sous-module de $\check{S}_p(E/D_n)$ engendré par les traces de $K[p^n]$ à $D_n$ des points de Heegner de niveau divisant $p^{n+1}$. On a alors la proposition :

5.6.5. **Proposition.** *La norme de $D_{n+1}$ à $D_n$ induit une application de $\mathcal{H}_{n+1}$ à $\mathcal{H}_n$. Elle est surjective pour $n \geq 1$. L'indice de $Tr_{n,0}(\mathcal{H}_n)$ dans $\mathcal{H}_0$ est égal à $L(E/K_p, 1)^{-1}$ (facteur d'Euler local en $p$).*

5.6.6. **Remarque.** La définition couramment admise est de prendre le sous-module de $\check{S}_p(E/D_n)$ engendré par les traces de $K[p^n]$ à $D_n$ des points de Heegner de niveau $p^{n+1}$. Malheureusement, ce n'est pas toujours gros ! L'énoncé de Mazur dans [20] est incorrect : la condition $a_p \equiv 2 \mod p$ est inutile avec cette définition et la surjectivité affirmée est fausse pour $a_p \equiv 1 \mod p$. Essentiellement, on a besoin des points de niveau $p$ et de niveau $1$ à la fois car la trace de $H[p]$ à $K$ de $y_p$ est un "multiple rationnel" (et non entier) de la trace de $H[1]$ à $K$.

5.6.2. *Cas où le rang de $E(K)$ est plus grand que 1.* Plaçons dans le cas où l'hypothèse de Heegner est vérifiée et où le rang de $E(K)$ est strictement supérieur à 1. Dans ce cas, on sait que l'image de $\mathcal{H}_\infty$ dans $E(K)$ est nul. Cependant, on peut cependant construire un élément de de $\mathbb{Z}_p \otimes E(K)$ de la manière suivante (à condition que $\text{III}(E/K)(p)$ soit fini, construction de Kolyvagin-Solomon) : on choisit un générateur $\gamma$ de $\text{Gal}(D_\infty/K)$ et $\gamma_n$ sa restriction à $D_n$. Soit un élément $z_\infty = (z_n)$ de $\mathcal{H}_\infty$ dont la projection est nulle (les $z_n$ se calculent en fonction des points de Heegner de niveau une puissance de $p$). Alors

$$\sum_{i=0}^{p^n-1} i\gamma_n^i z_n$$

converge dans $\varinjlim_n \check{S}(D_n)$ vers un élément de $\check{S}(K)$. Moins explicitement, cela revient à résoudre l'équation $z_\infty = (\gamma - 1)z'_\infty$ dans $\check{S}(D_\infty)$. Le fait que cette équation admet une solution vient de ce que vient de ce que l'application trace $\mathbb{Q}_p \otimes \check{S}(D_\infty)_{\text{Gal}(D_\infty/K)} \to \mathbb{Q}_p \otimes \check{S}_p(K)$ est un isomorphisme. Plus généralement, il est possible que l'image de $z'_\infty$ dans $E(K)$ soit encore nulle. Il existe un entier $r$ tel que $z_\infty = (\gamma - 1)^r z_\infty^{(r)}$ et tel que la projection de $z_\infty^{(r)}$ dans $\check{S}_p(K)$ soit non nulle. Cette projection donne alors un point non trivial de $\check{S}_p(K)$.

Soit $\epsilon$ le signe de l'équation fonctionnelle de $E/\mathbb{Q}$.

5.6.7. **Lemme.** *Si $z_\infty = (\gamma - 1)^r z_\infty^{(r)}$ et si $z_0^{(r)}$ est la projection de $z_\infty^{(r)}$ sur $E(K)$, on a $c(z_0^{(r)}) = -(-1)^r \epsilon z_0^{(r)} \mod \text{torsion}$.*

Cela se déduit de la relation $x^c = -\epsilon c \mod \text{torsion}$ pour un point de Heegner et du fait que $c\gamma c^{-1} = \gamma^{-1}$.

Prenons $r$ maximal. Si le rang de $E(\mathbb{Q})$ est impair, le rang de $E(K)^-$ est pair. Il serait intéressant de relier les parités de $r$ et de $\text{rg} \, E(\mathbb{Q})$ et même les valeurs. En tout cas, $z'_0$ appartient à $E(K)^{-(-1)^r \epsilon}$

## Références

[1] P. Billot. Quelques aspects de la descente sur une courbe elliptique dans le cas de réduction supersingulière. *Compositio Math.* 58 (1986), 341-369.

[2] M. I. Bashmakov, The cohomology of abelian varieties over a number field, *Russian Math. Survey* 27 (1972), 25-70.

[3] M. Bertolini. Selmer groups and Heegner points in anticyclotomic $\mathbb{Z}_p$-extensions. *Compositio Math.* 99 (1995), 153-182.

[4] M. Bertolini et H. Darmon. Kolyvagin's descent and Mordell-Weil groups over ring class fields. *J. reine angew. Math.* 412 (1990), 63-74.

[5] M. Bertolini et H. Darmon. Heegner points on Mumford-Tate curves. *Invent. Math.* 126 (1996), 413-456.

[6] M. Bertolini et H. Darmon. Iwasawa's main conjecture for elliptic curves over anticyclotomic $\mathbb{Z}_p$-extensions. Prépublication 2001

[7] Daniel Bertrand. : Sous-groupes à un paramètre $p$-adique de variétés de groupe. *Invent. Math.* 40, 1977, 171-193.

[8] C. Cornut. Mazur's conjecture on higher Heegner points. Prépublication 2001.

[9] M. Flach. A generalisation of the Cassels-Tate pairing J. Reine Angew. Math., 412 (1990), 113-127

[10] J. W. S. Cassels. Arithmetic on curves of genus 1. IV. Proof of the Hauptvermutung. *J. Reine Angew. Math.* 211 (1962), 95-112

[11] J.-M. Fontaine et B. Perrin-Riou. Autour des conjectures de Bloch et Kato : cohomologie galoisienne et valeurs de fonctions $L$ In *Motives* (Seattle, WA, 1991), Amer. Math. Soc. , Providence, RI (1994), 599–706.

[12] R. Greenberg. Iwasawa theory for elliptic curves. In *Arithmetic theory of elliptic curves* (Cetraro, 1997), p. 51-144. Springer, Berlin, 1999.

[13] R. Greenberg. Introduction to Iwasawa theory for elliptic curves. In *Arithmetic algebraic geometry* (Park City, UT, 1999) Amer. Math. Soc.Providence s(2001), 407-464.

[14] L. Guo. On a generalisation of Tate dualities with applications to Iwasawa theory. *Compos. Math.* 85 (1993), 125-161.

[15] Hida. A $p$-adic measure attached to the zeta function associated with two elliptic modular forms I. *Invent. Math.* 79 (1985), 159-195.

[16] U. Jannsen Iwasawa modules up to isomorphism. In *Algebraic number theory*, Adv. Stud. Pure Math., 17, Academic Press, Boston,pp. 171–207,

[17] V. A. Kolyvagin. Euler systems. In *The Grthendieck Festschrift*, Vol. II, Progress in Math. 87 (1990), Birkhäuser Boston, 435-483.

[18] H. Kato $p$-adic Hodge theory and values of zeta functions of modular forms prépublication (2001).

[19] W. McCallum A duality theorem in the multivariable Iwasawa theory of local fields. *J. Reine Angew. Math.* 464 (1995), 143-172.

[20] B. Mazur. Modular curves and arithmetic. in *Proceedings of the ICM* 1983 (Varsaw), Vol. 1, 185-211, PWN, Warsaw, 1984.

[21] J. Nekovář. On the parity of ranks of Selmer groups. II. *C. R. Acad. Sci. Paris Sér. I Math.* 332 (2001), no. 2, 99–104.

[22] J. Nekovář. Selmer complexes. preprint (2001)

[23] B. Perrin-Riou. Arithmétique des courbes elliptiques et théorie d'Iwasawa. Mém. Soc. Math. France (N.S.) No. 17 (1984), 130 pp

[24] B. Perrin-Riou. Fonctions L $p$-adiques associées à une forme modulaire et à un corps quadratique imaginaire, *J. London Math. Soc. (2)* 38 (1988), 1-32.

[25] B. Perrin-Riou. Fonctions $L$ $p$-adiques, théorie d'Iwasawa et points de Heegner. *Bull. Soc Math. France* 115 (1987), 399-456.

[26] B. Perrin-Riou. Théorie d'Iwasawa et hauteurs $p$-adiques. *Invent. Math.* 109 (1992), 137-185.

[27] D. E. Rohrlich. On $L$-functions of elliptic curves and anticyclotomic towers. *Invent. Math.* 75 (1984), 383-408

[28] V. Vatsal. Uniform distribution of Heegner points. *Invent. Math.* (2002)

[29] S-W Zhang Gross-Zagier Formula for $GL_2$ The Asian Journal of Mathematics, 5 (2001), 183-290.

Mathématiques, Bât. 425, Université Paris-Sud F-91405 Orsay, France
*E-mail address*: Bernadette.Perrin-Riou@math.u-psud.fr